\documentclass[12pt]{amsart}
\usepackage{amsmath,amsfonts,amsthm,amscd,}
\newtheorem{theorem}{Theorem}

\newtheorem{proposition}[theorem]{Proposition}
\newtheorem{lemma}[theorem]{Lemma}
\newtheorem{definition}[theorem]{Definition}

\newtheorem{corollary}[theorem]{Corollary}
\theoremstyle{remark}
\newtheorem{remark}[theorem]{Remark}
\newcommand{\ds}{\displaystyle}

\numberwithin{equation}{section} \numberwithin{theorem}{section}
\begin{document}
\bibliographystyle{amsalpha}
\title[Convexity]{Convexity and singularities of curvature equations
in conformal geometry}
\author{Matthew J. Gursky}
\address{Matthew J. Gursky\\
 Department of Mathematics \\ University of Notre Dame \\ Notre Dame, IN 46556}
\email{mgursky@nd.edu}
\thanks{The research of the first author was
partially supported by NSF Grant DMS-0500538.}
\author{Jeff Viaclovsky}
\address{Jeff Viaclovsky, Department of Mathematics, MIT, Cambridge, MA 02139}
\email{jeffv@math.mit.edu}
\thanks{The research of the second author was partially
supported by NSF Grant DMS-0503506.}
\begin{abstract}
We define a generalization of convex functions, which we call
$\delta$-convex functions, and show they must satisfy interior
H\"older and $W^{1,p}$ estimates.
As an application, we consider solutions of a certain class of fully
nonlinear equations in conformal geometry with isolated
singularities, in the case of non-negative Ricci curvature. We prove
that such solutions either extend to a H\"older continuous function
across the singularity, or else have the same singular behavior as
the fundamental solution of the conformal Laplacian.
We also obtain various removable singularity
theorems for these equations.
\end{abstract}
\maketitle

\section{Introduction}

Our goal in this paper is to understand the behavior of solutions to
certain geometric PDEs with isolated point singularities. Since the
relevant equations are fully nonlinear we impose an additional
condition on our solutions, known as {\it admissibility}, which
guarantees that the resulting equations are elliptic.  Although the
solutions are defined on deleted neighborhoods of Riemannian
manifolds, in the course of our analysis we are naturally lead to
the study of functions locally defined on Euclidean space which
satisfy a certain convexity condition. This notion of convexity is
weaker than the more familiar one of $k$-convexity, and has the
additional advantage of being a {\it linear} condition.

Let us begin by recalling some basic definitions.  For $1 \leq k
\leq n$, denote by $\sigma_k : \mathbf{R}^n \rightarrow
\mathbf{R}$ the $k$-th elementary symmetric polynomial, and
$\Gamma_{\sigma_k} \subset \mathbf{R}^n$ the component of $\{x \in
\mathbf{R}^n | \sigma_k(x)
> 0 \}$ containing the positive cone $\{ x \in \mathbf{R}^n | x_1
> 0,..., x_n > 0 \}$. A function $u \in C^2(\Omega)$, where
$\Omega \subset \mathbf{R}^n$ is open, is {\em $k$-convex in
}$\Omega$ if
\begin{align} \label{kconvex}
D^2u \in \Gamma_{\sigma_k}
\end{align}
at each point in $\Omega$. This condition naturally arises in the
study of the {\it Hessian equations}
\begin{align} \label{HessEqn}
\sigma_k(D^2u) = f(x).
\end{align}
In particular, (\ref{HessEqn}) is elliptic provided $u$ is
$k$-convex; see \cite{Garding}, \cite{CNSIII}, and
\cite{TrudingerWangII}, for example.

To introduce our new notion of convexity we need to define a
family of nested cones $\Gamma_{\delta} \subset \mathbf{R}^n$:

\begin{definition}  Let $ \delta \in \mathbf{R}$.
The $n$-tuple $\lambda =
(\lambda_1,\dots,\lambda_n) \in \Gamma_{\delta}$ if and only if
\begin{align}
\label{Ricposcond} \lambda_i > -\delta \sum_{j=1}^n \lambda_j \quad
\forall \ 1 \leq i \leq n,
\end{align}
and we define $ \overline{\Gamma}_{\delta}$ to be the closure of $\Gamma_{\delta}$.
\end{definition}

Note that $\Gamma_{\delta_1} \subset \Gamma_{\delta_2}$ whenever
$\delta_1 < \delta_2$, and also that
if $\lambda \in \Gamma_{\delta}$ for $\delta > -1/n$, then
\begin{align} \label{tracepos}
\sum_{i=1}^n \lambda_i > 0.
\end{align}

\begin{definition}  Let $W \in C^2(\Omega)$, where $\Omega \subset
\mathbf{R}^n$ is open.  We say that $W$ is strictly $\delta$-convex in
$\Omega$ if
\begin{align} \label{deltacondef}
D^2W \in \Gamma_{\delta}
\end{align}
at each point in $\Omega$,
and  $W$ is $\delta$-convex in $\Omega$ if
\begin{align} \label{deltacondef2}
D^2W \in \overline{\Gamma}_{\delta}
\end{align}
at each point in $\Omega$.
\end{definition}

For $ k > n/2$, a $k$-convex function is $\delta$-convex for
$\delta = \frac{n-k}{n(k-1)}$.  This fact is implicit in the work of
Trudinger-Wang \cite{TrudingerWangII}, where they proved {\it a
priori} estimates for $k$-convex functions.  However, they did not
express it in these terms. Our main results about
$\delta$-convex functions are corresponding H\"older and
$W^{1,p}$ estimates:
\begin{theorem}
\label{unifhold} Let $0 \leq \delta < \frac{1}{n-2}$,
$\Omega \subset \mathbf{R}^n$ be open, and assume
$W \in C^2(\Omega)$ is $\delta$-convex. Then on any domain $\Omega' \subset
\subset \Omega$, $W$ satisfies
\begin{align}
\label{estwl1} \| W \|_{C^{\gamma}(\Omega')} \leq C_1 \int_{\Omega} |W|,
\end{align}
where
\begin{align}
\gamma = \frac{ 1 +  (2-n)\delta}{ 1 + \delta},
\end{align}
and $C_1$ depends only upon $\Omega', \Omega,$ and $\delta$.
Furthermore, for $ p < p_{\delta} = \frac{ n(1 + \delta)}{(n-1)\delta}$,
\begin{align}
\label{w1pnew}
\left\{ \int_{\Omega'} |\nabla W|^p  \right\}^{1/p}
\leq C_2 \int_{\Omega} |W|,
\end{align}
and $C_2$ depends only upon $\Omega', \Omega,$ and $p$.
\end{theorem}

\begin{remark}  Theorem \ref{unifhold} is
a generalization of the classical theorem stating that convex
functions are Lipschitz \cite[Chapter 6]{EvansGariepy}. For
$k$-convex functions, the constants simplify to  $\gamma = 2 - n/k$
and $p_{\delta} = \frac{nk}{n-k}$, and our estimates agree with the
estimates of Trudinger and Wang \cite{TrudingerWangII}. While the
proof of Theorem \ref{unifhold} is based on the work of
Trudinger-Wang, we emphasize that only the condition of
$\delta$-convexity is used, which is weaker than $k$-convexity.
\end{remark}

\subsection{Conformal Geometry}
We now turn to the fully nonlinear equations in conformal geometry
which motivated this work.  To begin, let $(M^n,g)$ be a smooth,
closed Riemannian manifold of dimension $n$.  We denote the Ricci
tensor of $g$ by by $Ric$ and the scalar curvature by $R$.  In
addition, the {\it Weyl-Schouten tensor} is defined by
\begin{align}
\label{WStensor} A = \frac{1}{(n-2)}\Big( Ric - \frac{1}{2(n-1)} R
g\Big).
\end{align}
$A$ is a symmetric $(0,2)$-tensor; using the Riemannian metric we
can associate the dual tensor of type $(1,1)$ denoted $g^{-1}A$. In
classical language, $g^{-1}A$ is obtained from $A$ by ``raising an
index.''  The tensor $g^{-1}A$ can also be viewed as a symmetric
linear transformation of the tangent space at each point; thus it
has $n$ real eigenvalues.

Let $g_u = e^{-2u}g$ be a conformal metric and let $A_u$ denote the
Weyl-Schouten tensor of $g_u$. In this paper we will consider
singular solutions of equations of the form
\begin{align} \label{Fform}
F(g_u^{-1}A_u) = f(x),
\end{align}
where $F$ is a real-valued function of $n$-variables,
$F(g_u^{-1}A_u)$ means $F$ applied to the eigenvalues of
$g_u^{-1}A_u$, and $f(x)$ is a given function.  Since $A_u$ is
related to $A$ by the formula
\begin{align} \label{Achange}
A_u = A + \nabla^2 u + du \otimes du - \frac{1}{2}| \nabla u|^2 g
\end{align}
(see \cite{Jeff2}), in general (\ref{Fform}) will be a fully
nonlinear equation of second order.

An example of particular importance is when $F = \sigma_k^{1/k}$:
\begin{align}
\label{sigmak1} \sigma_k^{1/k}(g_u^{-1}A_u) = f(x).
\end{align}
Following the conventions of our previous papers \cite{GVAM},
\cite{GVRic} we use $g$ (not $g_u$) to raise the index in $A_u$.
 That is, we interpret $A_u$ as a bilinear form on the tangent space
with inner product $g$ (instead of $g_u$), and understand
$\sigma_k(\cdot)$ to mean $\sigma_k$ applied to the eigenvalues of
$g^{-1}A_u$.  Using this convention, equation (\ref{sigmak1})
becomes
\begin{align} \label{sigmak}
\sigma_k^{1/k}(A_u) = f(x)e^{-2u},
\end{align}
or, by (\ref{Achange}),
\begin{align}
\label{hessk} \sigma_k^{1/k}(A + \nabla^2 u + du \otimes du -
\frac{1}{2}| \nabla u|^2 g) = f(x)e^{-2u}.
\end{align}
Note that when $k=1$, then $\sigma_1(g^{-1}A) = trace(A) =
\frac{1}{2(n-1)}R$.  Therefore, (\ref{hessk}) is the scalar
curvature equation.

As in the example of the Hessian equations (\ref{HessEqn}), given
an open set $\Omega \subset M^n$ and a solution $u \in
C^2(\Omega)$ of (\ref{hessk}), $u$ is an elliptic solution if the
eigenvalues of $A_u$ are in $\Gamma_{\sigma_k}$ at each point of
$\Omega$.

As we observed above, if $A_u \in \Gamma_{\sigma_k}$ then $A_u \in
\Gamma_{\delta}$ for some $\delta = \delta(k,n) > 0$. Given an
open set $\Omega \subset M^n$ and $u \in C^2(\Omega)$, if the
eigenvalues of  $A_u$ are in $\Gamma_{\delta}$ at each point of
$\Omega$, we then say that $u$ is {\em strictly
$\delta$-admissible } in $\Omega$; if the eigenvalues are in  the
closure $\overline{\Gamma}_{\delta}$, then we say that $u$ is {\em
$\delta$-admissible} in $\Omega$.

It turns out that $\delta$-admissibility has an important
geometric consequence:  If the eigenvalues of the Schouten tensor
$A_g$ are in $\Gamma_{\delta}$ at each point of $M^n$, then
(\ref{tracepos}) for $\delta > -1/n$ implies the scalar curvature
of $(M^n,g)$ is positive, while (\ref{Ricposcond}) for
$\delta < \frac{1}{n-2}$ implies
the Ricci curvature is positive.  In fact,
\begin{align} \label{Ricposintro}
Ric_g - [ 1 + (2-n) \delta] \sigma_1(A_g)g \geq 0.
\end{align}
This fact is crucial in our analysis.

When considering the more general equation (\ref{Fform}) we need
to impose various structural conditions on the function $F$ and
its domain.  Suppose
\begin{align}
\label{Fex} F : \Gamma \subset \mathbf{R}^n \rightarrow \mathbf{R}
\end{align}
with $F \in C^{\infty}(\Gamma) \cap C^{0}(\overline{\Gamma})$,
where $\Gamma \subset \mathbf{R}^n$ is an open, symmetric, convex
cone.  In addition, we assume

\vskip.1in ${\bf (i)}$ $F$ is symmetric, concave, and homogeneous of
degree one.

\vskip.1in  ${\bf (ii)}$ $F > 0$ in $\Gamma$, and $F = 0$ on
$\partial \Gamma$.

\vskip.1in  ${\bf (iii)}$  $F$ is {\em elliptic}:
$F_{\lambda_i}(\lambda) > 0$ for each $1 \leq i \leq n$, $\lambda
\in \Gamma$.

\vskip.1in  ${\bf (iv)}$  $\Gamma \supset \Gamma_{\sigma_n}$, and
there exists a constant $0 \leq \delta< \frac{1}{n-2}$ such that $\Gamma
\subset \Gamma_{\delta}$.

\vskip.1in
\noindent For $F$ satisfying ${ \bf (i)-(iv) }$, consider the equation
\begin{align} \label{hessF}
F(A_u) = f(x)e^{-2u}.
\end{align}
Given an open set $\Omega \subset M^n$ and a solution $u \in
C^2(\Omega)$ of (\ref{hessF}), $u$ is an elliptic solution if the
eigenvalues of $A_u$ are in $\Gamma$ at each point of $\Omega$.  We
then say that $u$ is {strictly \it $\Gamma$-admissible} (or just
{\it strictly admissible}). By ${ \bf (iv) }$, any strictly
$\Gamma$-admissible solution is strictly $\delta$-admissible. We
will also be interested in solutions of (\ref{hessF}) with $f(x)
\geq 0$ and $A_u \in \overline{\Gamma}$. In this case equation
(\ref{hessF}) may be degenerate elliptic; correspondingly we say
such solutions are {\it $\Gamma$-admissible} (or just { \it
admissible}), and therefore $\delta$-admissible.

Some examples of interest are \vskip.1in \noindent {\bf Example 1.}
Let
\begin{align} \label{hesskF}
F(A_u) = \sigma_k^{1/k}(A_u) = f(x)e^{-2u}
\end{align}
with $\Gamma = \Gamma_{\sigma_k}$, $k > n/2$. Since $k > n/2$, by
\cite{GVW} we find that the eigenvalues of $A_u$ satisfy
inequality (\ref{Ricposcond}) with
\begin{align*}
0 \leq \delta = \frac{ n-k}{n(k-1)} < \frac{1}{n-2}.
\end{align*}
\vskip.1in \noindent {\bf Example 2.} Let $1 \leq l < k$ and $k
> n/2$, and consider
\begin{align} \label{rats}
F(A_u) = \ds \Big( \frac{\sigma_k(A_u)}{\sigma_l(A_u)}
\Big)^{\frac{1}{k-l}} = f(x)e^{-2u}.
\end{align}
In this case we also take $\Gamma = \Gamma_{\sigma_k}$. \vskip.1in
\noindent {\bf Example 3.} For $\tau \leq 1$ let
\begin{align*}
A^{\tau} = \frac{1}{(n-2)} \Big( Ric - \frac{\tau}{2(n-1)} R g
\Big),
\end{align*}
and consider the equation
\begin{align}
\label{sigmakt} F(A_u) = \sigma_k^{1/k}(A^{\tau}_u) = f(x)e^{-2u}.
\end{align}
By (\ref{Achange}), this is equivalent to the fully nonlinear
equation
\begin{align*}
\sigma_k^{1/k} \Big(
 A^{\tau} + \nabla^2 u + \frac{1-\tau}{n-2}(\Delta u)g
+ du \otimes du - \frac{2-\tau}{2} |\nabla u|^2 g \Big) =
f(x)e^{-2u}.
\end{align*}
In the Appendix of \cite{GVRic} we showed that the results of
\cite{GVW} imply the existence of $\tau_0 = \tau_0(n,k)
> 0$ and $\delta_0 = \delta(k,n) < \frac{1}{n-2}$ so that if
$1 \geq \tau > \tau_0(n,k)$ and  $A_g^{\tau} \in
\Gamma_{\sigma_k}$ with $k>n/2$, then $A_g \in \Gamma_{\delta}$
with $\delta = \delta_0.$ When we refer to this
example, we will tacitly make the assumption that $k > n/2$ and
$1 \geq \tau > \tau_0(n,k)$.

\vskip.2in In this paper we study solutions of (\ref{hessF}) with
isolated point singularities.   Thus, we assume $u$ is a solution of
(\ref{hessF}) in $\Omega = B(O,r_0) \setminus \{O\}$, and attempt to
understand the behavior of $u(x)$ as $x \to O$.

\begin{theorem}  \label{thm4}
Let $u \in C^4(\Omega)$ be an admissible solution of (\ref{hessF})
in $\Omega = B(O,r_0) \setminus \{O\}$, with $f(x) \equiv 0$ near $O$.
\vskip.1in \noindent $(i)$ If
\begin{align} \label{badinf2}
\liminf_{x \to O} u(x) > -\infty,
\end{align}
then $u$ can be extended to a H\"older continuous function $u^{\ast}
\in C^{\gamma}(B(O,r_0))$, with
\begin{align} \label{gammacon}
\gamma  =\frac{ 1 + (2-n)\delta}{1 + \delta}.
\end{align}
\vskip.1in \noindent $(ii)$ If
\begin{align} \label{badinf3}
\liminf_{x \to O} u(x) = -\infty,
\end{align}
then there is a constant $C > 0$ such that for all $x \in \Omega$,
\begin{align} \label{Greensgrowth}
2 \log d(x) - C \leq u(x) \leq 2 \log d(x) + C.
\end{align}
\end{theorem}

\begin{remark} In view of Proposition \ref{prop1} below (which only
uses the sign of the scalar curvature), a singular solution is
necessarily bounded from above.
\end{remark}

\begin{remark}
 In case $u$ satisfies (\ref{hesskF}) or
(\ref{rats}) (in particular, when $\Gamma = \Gamma_{\sigma_k}$ with
$k > n/2$), then
\begin{align}
\gamma = 2 - n/k > 0.
\end{align}
If $u$ satisfies (\ref{sigmakt}) with $k > n/2$ and $\tau > \tau_0 =
\frac{2(n-k)}{n}$, then
\begin{align}
\label{gamtdef} \gamma = \frac{ (n-2)(2k-2n+n\tau)}{ n-2k+kn-n\tau}
> 0.
\end{align}
\end{remark}
\begin{remark} We have a stronger statement of Theorem \ref{thm4} in the
locally conformally flat case; see Theorem \ref{thm1'} below.
\end{remark}
\begin{remark}
It is instructive to compare Theorem \ref{thm4} $(i)$ with the
recent classification of radial solutions of the
$\sigma_k$-curvature equations carried out by Chang-Han-Yang
\cite{CHY}. When $k
> n/2$ they show the existence of a solution to (\ref{hesskF}) in
$C^{2-n/k}(S^n)$, but whose second derivative blows up at an
isolated point.  Therefore, Theorem \ref{thm4} $(i)$ is optimal.
\end{remark}
\begin{remark}
To provide some context for the conclusions of part $(ii)$ in
Theorem \ref{thm4}, consider the case of the sphere with the round
metric $(S^n,g_{0})$.  Fix a point $O \in S^n$, and let $G$ denote
the Green's function for the conformal Laplacian $L$ with pole at
$O$. $G$ satisfies
\begin{align*}
L G = 0 \quad \mbox{on  } S^n \setminus \{ O \}, \\
G(x) \sim \mbox{dist}(x,O)^{2-n}.
\end{align*}
The the conformal manifold $(S^n \setminus \{ O \}, G^{
\frac{4}{(n-2)}}g_{0})$ is actually isometric to Euclidean space;
consequently the function $u = \ds -\frac{2}{(n-2)}\log G$ satisfies
\begin{align*}
F(A_{0} + \nabla^2 u + du \otimes du - \frac{1}{2}| \nabla u|^2
g_{0}) = 0
\end{align*}
on $S^n \setminus \{ O \}$. Thus, (\ref{Greensgrowth})
says, in essence, that a solution of (\ref{hessF})
with isolated singularity must blow up at the
same rate as the fundamental solution of the conformal Laplacian.
In fact, using the results from Section 7 of \cite{GVRic}, it 
is possible to show that (\ref{Greensgrowth}) implies furthermore that 
the background metric is the Euclidean metric, and 
$u(x) \equiv 2 \log r + C$. We omit the details, but the main idea is to look at 
$(M^n \setminus \{ O \}, g_u)$ as a complete manifold with 
non-negative Ricci curvature outside of a compact 
set. An adaptation of the arguments in Section 7 of \cite{GVRic} 
imply that the end is asymptotically flat, and 
applying a version of Bishop's volume comparison theorem, we find 
that $(M^n \setminus \{ O \}, g_u)$ is isometric to Euclidean space
outside of a compact set.  
\end{remark}

For the lower bound in (\ref{Greensgrowth}) above, we can weaken the regularity
assumption to
$u \in C^3$, but we must specialize to the $F$ given
in the above examples:
\begin{theorem}  \label{thm4q}
Let $u \in C^3(\Omega)$ be an admissible solution 
of either (\ref{hesskF}), (\ref{rats}), or (\ref{sigmakt})
in $\Omega = B(O,r_0) \setminus \{O\}$, with $f(x) \equiv 0$ near $O$.

If
\begin{align} \label{badinf3q}
\liminf_{x \to O} u(x) = -\infty,
\end{align}
then there is a constant $C > 0$ such that for all $x \in \Omega$,
\begin{align} \label{Greensgrowthq}
2 \log d(x) - C \leq u(x).
\end{align}
\end{theorem}
\begin{remark}
The proof of the estimate in (\ref{Greensgrowthq}) requires a local gradient
estimate for $C^3$ solutions,
which is currently only known for the special
cases in Examples 1--3. We conjecture it is true for
general $F$. If $u \in C^4$, then a local $C^2$-estimate
for general $F$ follows from
Sophie Chen's work \cite{Sophie1}, which was used in Theorem \ref{thm4}.
This is discussed in detail in Section \ref{Pointwise}.
\end{remark}

\subsection{Relation to the Existence Theory}
Theorem \ref{thm4} is in fact a corollary of a much more general
result about the local behavior of admissible functions.  This
result depends upon an explicit but rather subtle relationship
between admissibility and $\delta$-convexity; see Section
\ref{adtocon}.  Moreover, this generalization of Theorem \ref{thm4}
holds under weaker regularity assumptions--an important
consideration for certain applications, for reasons we now explain.

Aside from its intrinsic interest, the study of solutions with
isolated singularities is central to the study of {\it a priori}
estimates for solutions of (\ref{hessF}), and the related problem of
analyzing the blow-up of sequences of solutions. Both of these
topics were treated in our previous paper \cite{GVRic}, where we
proved a general existence result for solutions of (\ref{hessF})
assuming properties ${ \bf (i)-(iv)}$ above and certain {\it a
priori} estimates are satisfied.

Precisely because singular solutions often appear as limits of
smooth ones, there is an additional technical difficulty that often
arises. Namely, the limit may only be in $C^{1,1}_{loc}$ and satisfy
(\ref{hessF}) almost everywhere. For example, in \cite{GVRic}, a
divergent sequence of solutions $\{ u_i \}$ to (\ref{hessF}) is
rescaled by defining $v_i = u_i + \tau_i$, where $\{ \tau_i \}$ is a
sequence of numbers with $\tau_i \to +\infty$ as $i \to \infty$.
Each $v_i$ is also a solution of (\ref{hessF}), but with $f_i(x) =
e^{-2\tau_i}f(x)$.  Now, the sequence $\{ v_i \}$ converges (away
from a finite point set $\Sigma$), but the limit $v \in
C^{1,1}_{loc}(M^n \setminus \Sigma)$ is a (degenerate) admissible
solution of (\ref{hessF}) with $f(x) \equiv 0$. Thus, when studying
singular solutions of (\ref{hessF}) it is natural to impose the
weakest possible regularity.

A similar construction, by the way, was carried out in Schoen's work
on the Yamabe problem \cite{Schoen1}.  In this case, a divergent
minimizing sequence for the Yamabe functional is rescaled, and a
subsequence converges (away from a finite point set $\Sigma$) to a
solution of
\begin{align} \label{green}
Lh = 0 \quad \mbox{on  } M^n \setminus \Sigma,
\end{align}
where $L = \Delta - \frac{(n-2)}{4(n-1)}R$.  The important
difference here is that while $h$ is a singular solution of
(\ref{green}), it is {\it smooth} away from the singular points.
This allows one to apply the results of Serrin \cite{Serrin} and
Gilbarg-Serrin \cite{GS}, who classified $C^2$-solutions of
(\ref{green}) with isolated singularities: in fact, $h$ must be a
linear combination of fundamental solutions of the conformal
laplacian.

In general, to understand the behavior of solutions near isolated
singularities some form of the Harnack inequality seems essential,
as it was in the work of Gilbarg-Serrin for the semilinear case.
While Harnack inequalities have been established for solutions of
(\ref{hessF}) (see \cite{Jeff2}, \cite{LiLi1}, \cite{GuanWang1},
\cite{Sophie1})
they all assume at the very least $u \in C^3$ and $f \in C^1$, for
the simple reason that the proofs rely on differentiating the
equation.  In our existence work \cite{GVRic} described above, we
were able to show that the singular solution $v$ satisfied a Harnack
inequality by using the fact it was the limit of smooth solutions,
each of which satisfied the local gradient bounds proved by Guan and
Wang \cite{GuanWang1}. But for an arbitrary solution of
(\ref{hessF}) in, say, $C_{loc}^{1,1}$, it remains an open question
whether one can obtain such an estimate.

\subsection{Scale-Invariant Estimates}
To clarify precisely what is lacking for weak solutions of
(\ref{hessF}), we introduce the following terminology:

\begin{definition}  \label{sige} Let $B = B(O,r_0) \subset M^n$ and $\Omega = B
\setminus \{O\}$.  For $x \in \Omega$, let $d(x)$ denote the
distance to $O$. We say that $u \in C^{1}_{loc}(\Omega)$ satisfies a
{\em scale-invariant $C^1$-estimate} if there is a constant $C$ such
that
\begin{align} \label{gradr1}
|\nabla u(x)| \leq \ds \frac{C}{d(x)},
\end{align}
for every $x \in \Omega$. We say that $u \in C^{1,1}_{loc}(\Omega)$
satisfies a {\em scale-invariant $C^2$-estimate} if there is a
constant $C$ such that
\begin{align} \label{gradr}
|\nabla u(x)|^2 + |\nabla^2 u(x)| \leq \ds \frac{C}{d(x)^2}
\end{align}
for almost every $x \in \Omega$.
\end{definition}

\begin{remark}  Suppose $u$ is a solution of
(\ref{hessF}) on $\mathbf{R}^n \setminus \{ 0 \}$ with $f(x) \equiv
0$.  If $u$ satisfies (\ref{gradr}), then so does $u_{\lambda}(x) =
u(\lambda x)$ for any $\lambda > 0$ (with the same constant $C$).
This is why (\ref{gradr}) is called a ``scale-invariant'' estimate.
For example, consider $u \in C^{\infty}(\mathbf{R}^n \setminus \{ 0
\})$ given by
\begin{align*}
u(x) = \log |x|^2.
\end{align*}
Then $u$ is a solution of (\ref{hessF}) with $f(x) \equiv 0$.
Moreover, $|\nabla u_{\lambda}(x)| = |\nabla u(\lambda x)|$.
\end{remark}

We will postpone for now the question of when a scale-invariant
estimate can be verified.  Instead, we will first restate our main
result for $C_{loc}^{1,1}$-solutions of (\ref{hessF}), with
(\ref{gradr1}) and (\ref{gradr}) as additional assumptions:

\begin{theorem} \label{thm1}
Suppose $u \in C^{1,1}_{loc}$ satisfies a scale-invariant
$C^2$-estimate, and $A_u \in \overline{\Gamma}_{\delta}$ almost
everywhere in $\Omega = B(O,r_0) \setminus \{O\}$, where $0 \leq
\delta < \frac{1}{n-2}$. Then the conclusions of Theorem \ref{thm4}
hold.

\end{theorem}

When the background metric is locally conformally flat, we only need
a scale-invariant $C^1$-estimate:

\begin{theorem} \label{thm1''}  Assume $(M,g)$ is locally
conformally flat.  Suppose $u \in C^{1,1}_{loc}$ satisfies a
scale-invariant $C^1$-estimate, and $A_u \in
\overline{\Gamma}_{\delta}$ almost everywhere in $\Omega = B(O,r_0)
\setminus \{O\}$, where $0 \leq \delta < \frac{1}{n-2}$. Then the
conclusions of Theorem \ref{thm4} hold.
\end{theorem}

In fact, part $(i)$ of Theorem \ref{thm4} holds without assuming a
scale-invariant estimate; for part $(ii)$, we can verify a one-sided
bound:

\begin{theorem} \label{thm1'}  Assume $(M,g)$ is locally
conformally flat.  Suppose $u \in C^{1,1}_{loc}$ with $A_u \in
\overline{\Gamma}_{\delta}$ almost everywhere in $\Omega = B(O,r_0)
\setminus \{O\}$, where $0 \leq \delta < \frac{1}{n-2}.$ \vskip.1in
\noindent $(i)$ If
\begin{align} \label{goodinf'}
\liminf_{x \to O} u >  -\infty,
\end{align}
then $u$ can be extended to a H\"older continuous function $u^{\ast}
\in C^{1,1}_{loc}(\Omega) \cap C^{\gamma}(B(O,r_0))$,
 with $\gamma$ given by (\ref{gammacon}).

\vskip.1in \noindent  $(ii)$  If
\begin{align} \label{badinf'}
\liminf_{x \to O} u(x) = -\infty,
\end{align}
then there is a constant $C > 0$ such that for all $x \in \Omega$,
\begin{align} \label{neargrowth'}
u(x) \leq 2 \log d(x) + C.
\end{align}
\end{theorem}

\begin{remark}
The main reason we have a stronger statement in the locally
conformally flat case is roughly that, in normal coordinates, a
general Riemannian metric will be close to Euclidean only to second
order, while in the locally conformally flat case, we can find a
conformal metric which is exactly Euclidean in a neighborhood of a
point. More precisely, if we let $\{ x^i \}$ denote normal coordinates
(with respect to the background metric $g$) centered at $O$,
the cone condition  $A_u \in \overline{\Gamma}_{\delta}$
is equivalent to
\begin{align} \label{localcone3intro}
g^{jl}\left( \partial_i \partial_j u - \Gamma_{ij}^k \partial_k u +
u_i u_j - \frac{1}{2}g^{pq} \partial_p u \partial_q u g_{ij} +
A_{ij} \right) \in \overline{\Gamma}_{\delta}
\end{align}
a.e. in $\Omega$. In normal coordinates, $g^{jl} = \delta^{jl} +
O(|x|^{2}),$ and $|\Gamma_{ij}^k| = O(|x|)$ as $|x| \rightarrow 0$.
In particular, we find error terms of the form
$ O(|x|^2) ( \partial_i \partial_j u)$ as $|x| \rightarrow 0$,
which could be unbounded without a scale-invariant $C^2$-estimate on $u$.
\end{remark}

\vskip.2in

We now turn to the question: when does a solution of (\ref{hessF})
satisfy a scale-invariant $C^2$-estimate? If $u \in C^4(\Omega)$ is
a solution of (\ref{hessF}), then the local estimates proved by
Sophie Chen \cite{Sophie1} can be used to verify a scale-invariant
$C^2$-estimate.  If $u \in C^3(\Omega)$ and is a solution of either
(\ref{hesskF}), (\ref{rats}), or (\ref{sigmakt}), then the local
estimates of solutions established in \cite{GuanWang1},
\cite{GuanWang2}, \cite{LiLi2}, and \cite{GVJDG} can be used to to
verify the scale-invariant $C^1$-estimate:

\begin{proposition}
\label{prop3}
\noindent $(i)$ Let $u \in C^4(\Omega)$ be an
admissible solution of (\ref{hessF}) in  $\Omega = B(O,r_0)
\setminus \{O\}$. If $f \equiv 0$ in a neighborhood of $O$, then $u$
satisfies the scale-invariant $C^2$-estimate (\ref{gradr}).

\vskip.1in \noindent $(ii)$  Let $u \in C^3(\Omega)$ be an
admissible solution of either (\ref{hesskF}), (\ref{rats}), or
(\ref{sigmakt}) in $\Omega = B(O,r_0) \setminus \{O\}$.  If $f
\equiv 0$ in a neighborhood of $O$, then $u$ satisfies the
scale-invariant $C^1$-estimate (\ref{gradr1}).
\end{proposition}

\begin{remark}  In particular, Theorem \ref{thm4} follows from Proposition \ref{prop3} and Theorem
\ref{thm1}.  \end{remark}

\subsection{H\"older extension}

A geometrically natural condition to consider is that of finite
volume.  For example, suppose $u \in C^{1,1}_{loc}$ satisfies the
hypotheses of Theorem \ref{thm1}, and the volume of $g_u = e^{-2u}g$
is finite:
\begin{align} \label{volbound}
Vol_{g_u}(\Omega) = \int_{\Omega} e^{-nu} dvol_g < \infty.
\end{align}
Then by examining the integrand in (\ref{volbound}), it is clear
that $u$ cannot satisfy (\ref{Greensgrowth}).  Consequently, we have

\begin{corollary} \label{cor1}   Let $0 \leq \delta < \frac{1}{n-2}$,
$u \in C^{1,1}_{loc}$ satisfy $A_u \in \overline{\Gamma}_{\delta}$
almost everywhere in $\Omega = B(O,r_0) \setminus \{ O \}.$  Assume
that $g$ is locally conformally flat and $u$ satisfies a
scale-invariant $C^1$-estimate, or that $u$ satisfies a
scale-invariant $C^2$-estimate in the general case. If the volume of
the conformal metric $g_u = e^{-2u}g$ is finite, then $u$ can be
extended to a H\"older continuous function $u^{\ast} \in
C^{1,1}_{loc}(\Omega) \cap C^{\gamma}(B(O,r_0))$, where $\gamma$ is
given in  (\ref{gammacon}).
\end{corollary}

In addition, finite volume actually implies a scale-invariant
estimate:

\begin{theorem} \label{thm5}
Let $u \in C^4$. Assume that $u$ is an admissible solution of
(\ref{hessF}) in $\Omega = B(O,r_0) \setminus \{O\}$. If the volume
of the conformal metric $g_u = e^{-2u}g$ is finite; i.e., if
(\ref{volbound}) is satisfied, then $u$ satisfies a scale-invariant
$C^2$-estimate. Consequently, by Corollary \ref{cor1}, $u$ can be
extended to a H\"older continuous function $u^{\ast} \in
C^{\gamma}(B(O,r_0))$.
\end{theorem}

In the locally conformally flat case, when $u$ is a solution of one
the special examples (\ref{hesskF})--(\ref{sigmakt}), then we can
slightly weaken the regularity assumption:

\begin{theorem} \label{thm5'}
Let $(M^n,g)$ be locally conformally flat and $u \in C^{3}$. Assume
that $u$ is an admissible solution of either (\ref{hesskF}),
(\ref{rats}), or (\ref{sigmakt}) in $\Omega = B(O,r_0) \setminus
\{O\}$. If the volume of the conformal metric $g_u = e^{-2u}g$ is
finite then $u$ satisfies a scale-invariant $C^1$-estimate.
Consequently, by Corollary \ref{cor1}, $u$ can be extended to a
H\"older continuous function $u^{\ast} \in C^{\gamma}(B(O,r_0))$.
\end{theorem}

In \cite{Gonzalez2} Gonzalez studied the behavior of solutions to
(\ref{hessk}), $k < n/2$, with isolated singularities.  She proved
that $C^3$-solutions with finite volume are bounded across the
singularity.  In related work, Han \cite{Han} proved local
$L^{\infty}$-estimates for $W^{2,2}$-solutions of (\ref{hessF}) when
$k = 2$ and $n = 4$, assuming a smallness condition on the volume.

In a subsequent paper \cite{Gonzalez3} Gonzalez considered a
subcritical version of (\ref{hessk}):
\begin{align} \label{subhessk}
\sigma_k^{1/k}(A_u) = f(x)e^{-2 \beta u},
\end{align}
where $\beta < 1$.  Solutions of (\ref{subhessk}) with isolated
singularities are either bounded or satisfy a (sharp) growth
condition near the singularity analogous to (\ref{Greensgrowth}).

In the conformally flat case, when $f(x) \geq c_0 > 0$ near the
singularity $O$ we can also rule out blow-up:

\begin{theorem} \label{thm6}
Let $(M^n,g)$ be locally conformally flat and $u \in C^3$ a
strictly admissible solution 
of either (\ref{hesskF}), (\ref{rats}), or (\ref{sigmakt})
in $\Omega = B(O,r_0) \setminus
\{O\}$, with $f(x) \geq c_0 > 0$ near $O$.  Then $u$ can be
extended to a H\"older continuous function $u^{\ast} \in
C^3(\Omega) \cap C^{\gamma}(B(O,r_0))$.
\end{theorem}
\vskip.1in

\begin{remark}
This is proved in Section \ref{fgeq0}.
We conjecture that Theorem \ref{thm6} is true without the local
conformal flatness assumption.
\end{remark}

\subsection{H\"older and $L^p$-estimates}

Our final results are some $W^{1,p}-$ and H\"older-estimates for
$\delta$-admissible functions:

\begin{theorem}
\label{holderthm} Let $0 \leq \delta < \frac{1}{n-2}$, $u \in
C^{1,1}_{loc}$ satisfy $A_u \in \overline{\Gamma}_{\delta}$ almost
everywhere in $B(O,r_0).$ Let $v = e^{\beta u}$, where
\begin{align}
\label{betadef} \beta = \frac{ 1 + (2-n)\delta}{2(1 + \delta)}.
\end{align}
Then for $ p < p_{\delta} = \frac{ n(1 + \delta)}{(n-1)\delta}$,
\begin{align}
\label{optab} \left\{ \int_{B(x_0,r/2)} |\nabla v|^p  \right\}^{1/p}
\leq C \int_{B(x_0,r)} |v|,
\end{align}
where $C = C(r, p, \Vert g \Vert_{C^2})$. Consequently, if $ \gamma
= \frac{ 1 + (2-n)\delta}{1 + \delta}$, then by the Sobolev
embedding theorem
\begin{align}
\label{hest2first} \| v \|_{C^{\alpha}(B(x_0,r/2))} \leq C
\int_{B(x_0,r)} |v|,
\end{align}
for any $\alpha < \gamma$, where $C = C( r, \alpha, n , \Vert g
\Vert_{C^2})$. If $g$ is locally conformally flat, then we may take
$\alpha = \gamma$ in (\ref{hest2first}).
\end{theorem}

 We mention a related body of work which considers
solutions of (\ref{hessk}) defined on subdomains in the sphere.  By
the work of Schoen and Yau \cite{SchoenYau}, such solutions arise
when considering complete, conformally flat (admissible) metrics.
The goal is to study the singular set and derive estimates for the
Hausdorff dimension (see \cite{ChangHangYang1}, \cite{Gonzalez1},
\cite{GuanLinWang3}).

 In closing, we would also like to mention the preprints of
Yanyan Li \cite{Yanyana}, \cite{Yanyanb}, and Trudinger-Wang \cite{TWnew} 
which contain closely related results regarding the above H\"older
estimates, and asymptotic behaviour of solutions at singularities.

\subsection{Acknowledgements}
The authors would like to thank Alice Chang, Pengfei Guan, Yanyan Li,
Guofang Wang, and Paul Yang for helpful discussions and valuable remarks.
We are very grateful to Sophie Chen for informing us
of her nice work \cite{Sophie1}. 

\section{H\"older estimates for $\delta$-convex functions}
\label{holsec}

We begin by giving the proof of the H\"older estimate (\ref{estwl1})
in Theorem \ref{unifhold}.
The $W^{1,p}$ estimate (\ref{w1pnew}) will be proved later in
Section \ref{inthol}. First, notice that the function
\begin{align}
G(x) = C |x-y|^{\gamma}
\end{align}
is a solution of
\begin{align}
\label{gammq} F_{\delta}( D^2 G) \equiv \mbox{det}^{1/n} \Big( D^2 G
+ \delta (\Delta G) I \Big) = 0
\end{align}
on $R^n - \{ y \}$, with $D^2 G \in \overline{\Gamma}_{\sigma_n}$,
for
\begin{align}
\gamma = \frac{ 1 +  (2-n)\delta}{ 1 + \delta}.
\end{align}
To see this, by scaling and translation, assume that $C=1$, and $y =
0$, then
\begin{align}
\partial_i \partial_j G &= \gamma( \gamma-2) |x|^{\gamma-4}x_i x_j
+ \gamma |x|^{\gamma-2} \delta_{ij},\\
\Delta G &= ( \gamma ( \gamma-2) + n \gamma) |x|^{\gamma-2}.
\end{align}
A computation shows that
\begin{align*}
\partial_i \partial_j G
+ \delta (\Delta G) \delta_{ij}&=
\gamma( \gamma-2) |x|^{\gamma-4} x_i x_j \\
& \ \ \ \ \ \ + \Big( \gamma + \delta ( \gamma (
\gamma-2) + n \gamma) \Big)
|x|^{\gamma-2} \delta_{ij},\\
& = \gamma( \gamma-2) |x|^{\gamma-4} x_i x_j - \gamma( \gamma-2)
|x|^{\gamma-2} \delta_{ij}.
\end{align*}
This has one zero eigenvalue, and the other eigenvalues are all
equal to $ \gamma( 2 - \gamma) > 0$, so $(\ref{gammq})$ follows,
with  $D^2 G \in \overline{\Gamma}_{\sigma_n}$.

Equation (\ref{gammq}) is a fully nonlinear (degenerate) elliptic
equation, with concave $F_{\delta}$ and ellipticity cone
$\overline{\Gamma}_{\sigma_n}$. Let $ y \in B$, $r>0$, and define on
$B(y, R) \setminus \{ y \}$,
\begin{align}
\tilde{W}(x) &=
\frac{W(x) - W(y)}{osc_{B_R(y)}W },\\
G(x) &=  \left( \frac{|x-y|}{R} \right)^{\gamma}.
\end{align}
Choose $ \epsilon > 0$; we have
\begin{align}
\label{strell} F_{\delta}(D^2 (\tilde{W}(x) - \epsilon )) &> 0
\mbox{ on }
B(y, R) \setminus B(y,r), \\
\label{degell}
F_{\delta}(D^2 G) & = 0 \mbox{ on } B(y, R) \setminus B(y,r),\\
\tilde{W} - \epsilon &\leq G \mbox{ on }
\partial (B(y, R) \setminus B(y,r)),
\end{align}
for $r$ sufficiently small. Note that (\ref{strell}) is strictly
elliptic, but (\ref{degell}) is degenerate elliptic. By Theorem
\cite[17.1]{GT}, the difference $\tilde{W}(x) - \epsilon  - G$
satisfies a linear strictly elliptic equation. From the maximum
principle,
\begin{align}
\tilde{W} -\epsilon \leq G \mbox{ on } B(y, R) \setminus B(y,r).
\end{align}
Since $\epsilon > 0$ is arbitrary, for $x \in B_R(y) \setminus \{y
\}$ we obtain
\begin{align}
W(x) - W(y) \leq (osc_{B_R(y)} W) \left( \frac{|x-y|}{R}
\right)^{\gamma}.
\end{align}
The estimate (\ref{estwl1}) then follows by a standard interpolation
argument; see \cite{TrudingerWangII}.

\section{Admissibility and $\delta$-convexity} \label{adtoco}

The next result explains the relationship between admissibility
and $\delta$-convexity:

\begin{theorem} \label{adtocon} Suppose $A_u \in \overline{\Gamma}_{\delta}$,
where $A_u$ is given by (\ref{Achange}).
Let $v = e^{\beta u}$, where
\begin{align}
\beta = \frac{ 1 + (2-n)\delta}{2(1 + \delta)}.
\end{align}
Then
\begin{align} \label{hessvcone}
\nabla^2 v + \beta v A_g \in \overline{\Gamma}_{\delta}.
\end{align}
In particular, if $g_u = e^{-2u}ds^2$, where $ds^2$ is the flat metric on $\mathbf{R}^n$, then $v$ is $\delta$-convex.
\end{theorem}

\begin{proof}  Since $\log v = \beta  u$, we have
\begin{align}
\beta \nabla u &= \frac{\nabla v}{v},\\
\beta \nabla^2 u &= \frac{\nabla^2 v }{v} - \frac{dv \otimes
dv}{v^2}.
\end{align}
Letting $\alpha = \beta^{-1}$ and using (\ref{Achange}), we obtain
\begin{align}
\label{change2b} {A}_v &= A + \alpha \frac{ \nabla^2 v }{v} + (
\alpha^2 - \alpha) \frac{1}{v^2} dv \otimes dv - \frac{1}{2}
\alpha^2 \frac{ |\nabla v|^2}{v^2} g.
\end{align}

In terms of $v$, the admissibility condition $A_v \in
\overline{\Gamma}_{\delta}$ implies
\begin{align} \label{vcone}
v A_g + \alpha v \nabla^2 v + ( \alpha^2 - \alpha) dv \otimes dv -
\frac{\alpha^2}{2}|\nabla v|^2g \in \overline{\Gamma}_{\delta}.
\end{align}
Now examine the gradient terms, which are proportional to:
\begin{align*}
 ( \alpha - 1) dv \otimes dv - \frac{\alpha}{2}|\nabla v|^2g.
\end{align*}
In terms of $\delta$, this is
\begin{align*}
 \frac{1 + n \delta}{1 + (2-n)\delta} dv \otimes dv -
\frac{ 1 + \delta}{1 + (2-n)\delta}|\nabla v|^2g,
\end{align*}
which is proportional to
\begin{align}
\label{gradthing}
 \frac{ 1 + n \delta}{1 + \delta} dv \otimes
dv - |\nabla v|^2g.
\end{align}
The eigenvalues of this tensor are
\begin{align}
\left( \frac{ (n-1) \delta}{1 + \delta}, -1, \dots, -1
\right)|\nabla v|^2,
\end{align}
and the trace is
\begin{align*}
 -\frac{n-1}{1 + \delta} |\nabla v|^2.
\end{align*}
Clearly, this implies that (\ref{gradthing}) belongs to $-
\overline{\Gamma}_{\delta}$. Since $\overline{\Gamma}_{\delta}$ is a
convex cone, it follows that
\begin{align} \label{hessvconegoog}
\nabla^2 v + \beta v A_g \in \overline{\Gamma}_{\delta}.
\end{align}
\end{proof}

\section{Preliminary estimates for singular solutions}

In this section we prove some technical results which will be used in the proofs of the main
theorems.

\subsection{Pointwise estimates}
The following Proposition gives an upper bound for solutions
in $B(O,r_0) \setminus \{O\}$, assuming
only that the scalar curvature is non-negative.
\begin{proposition} \label{prop1} Let $u \in C^{1,1}_{loc}(\Omega)$
satisfy
\begin{align} \label{scapos}
\sigma_1(A + \nabla^2 u + du \otimes du - \frac{1}{2}| \nabla u|^2
g) \geq 0
\end{align}
almost everywhere in $\Omega = B(O,r_0) \setminus \{O\}$.  Then
\begin{align} \label{supbound}
\sup_{\Omega} u < +\infty.
\end{align}
\end{proposition}
\begin{proof}
The inequality (\ref{scapos}) is simply
\begin{align*}
\Delta u - \ds \frac{(n-2)}{2}|\nabla u|^2 + \sigma_1(A) \geq 0.
\end{align*}
a.e. in $\Omega$, or
\begin{align} \label{sce}
\Delta u \geq \ds \frac{(n-2)}{2}|\nabla u|^2 - \sigma_1(A)
\end{align}
a.e. in $\Omega$.

Our first observation is

\begin{lemma} \label{uinw12}  Let $B = B(O,r_0)$.  Then $u \in
W^{1,2}(B)$.
\end{lemma}

\begin{proof}  For $\epsilon > 0$ small, let $\eta_{\epsilon}$ denote a cut-off
function supported in $\Omega$ satisfying
\begin{align} \label{etadef}
\eta_{\epsilon}(x) = \left\{ \begin{array}{lll}
0 & x \in B(O,\epsilon), \\
1 & x \in B(O,r_0/4) \setminus B(O,2\epsilon), \\
0 & x \in B \setminus B(O,r_0/2),
\end{array}
\right.
\end{align}
and $|\nabla \eta_{\epsilon}| \leq C/\epsilon.$  Since (\ref{sce})
holds a.e. on $\Omega$ and $\eta_{\epsilon}^2$ is supported in
$\Omega$ we have
\begin{align} \label{weak1}
\int \eta_{\epsilon}^2 \Delta u \geq \int \ds \frac{(n-2)}{2}|\nabla
u|^2 \eta_{\epsilon}^2 - \int \eta_{\epsilon}^2 \sigma_1(A).
\end{align}
Integrating by parts,
\begin{align*}
-\int 2 \langle \nabla \eta_{\epsilon}, \nabla u \rangle
\eta_{\epsilon}  \geq \int \ds \frac{(n-2)}{2}|\nabla u|^2
\eta_{\epsilon}^2  - \int \eta_{\epsilon}^2 \sigma_1(A).
\end{align*}
Using the inequality
\begin{align*}
-\int 2 \langle \nabla \eta_{\epsilon}, \nabla u \rangle
\eta_{\epsilon} \leq \int \ds \frac{(n-2)}{4}|\nabla u|^2
\eta_{\epsilon}^2 + \int \ds \frac{4}{(n-2)} |\nabla
\eta_{\epsilon}|^2,
\end{align*}
we conclude
\begin{align*}
\int \ds \frac{(n-2)}{2}|\nabla u|^2 \eta_{\epsilon}^2  - \int
\eta_{\epsilon}^2 \sigma_1(A) \leq \int \ds \frac{(n-2)}{4}|\nabla
u|^2 \eta_{\epsilon}^2 + \int \ds \frac{4}{(n-2)} |\nabla
\eta_{\epsilon}|^2,
\end{align*}
which implies
\begin{align} \label{dubd}
\int |\nabla u|^2 \eta_{\epsilon}^2 \leq C\int \Big[ |\nabla
\eta_{\epsilon}|^2 + \eta_{\epsilon}^2 \Big].
\end{align}
Note that
\begin{align*}
\int |\nabla \eta_{\epsilon}|^2 &\leq C \epsilon^{-2}
\int_{B(O,2\epsilon) \setminus B(O,\epsilon)}  \\
&\leq C\epsilon^{(n-2)} \to 0
\end{align*}
as $\epsilon \to 0$.  Therefore, letting $\epsilon \to 0$ in
(\ref{dubd}), we get
\begin{align} \label{dul2}
\int_{B} |\nabla u|^2 \leq C.
\end{align}

To prove that $u \in L^2(B)$ we apply the Poincare inequality, which
states
\begin{align}   \label{pc}
\int_{B} \varphi^2 \leq \ds \frac{1}{\lambda_1} \int_{B} |\nabla
\varphi|^2
\end{align}
for all $\varphi \in W_0^{1,2}(B)$, where $\lambda_1$ is the first
(Dirichlet) eigenvalue of $-\Delta$ on $B$. For $k \geq 1$, define
\begin{align} \label{ukdef}
u_k(x) = \left\{ \begin{array}{lll}
k & \mbox{if}\quad u(x) \geq k, \\
u(x) & \mbox{if}\quad k \geq u(x) \geq -k, \\
-k & \mbox{if}\quad u(x) \leq -k.
\end{array}
\right.
\end{align}
Let $\zeta$ be another cut-off function supported in $B$, this time
satisfying
\begin{align} \label{zetaprop}
\zeta(x) \equiv 1 \quad \forall x \in B(O,r_0/4).
\end{align}
For each $k \geq 1$, (\ref{dul2}) implies $|\nabla u_k| \in L^2(B)$,
and since $u_k$ is bounded, it follows that $u_k \in W^{1,2}(B)$.
Therefore, by (\ref{pc}),
\begin{align} \label{pic1} \begin{split}
\int_{B} \zeta^2 u_k^2 &\leq \ds \frac{1}{\lambda_1}
\int_{B} |\nabla (\zeta u_k)|^2 \\
&\leq 2 \Big[ \int_{B} \zeta^2 |\nabla u_k|^2 + \int_{B} u_k^2
|\nabla \zeta|^2 \Big]. \end{split}
\end{align}
By (\ref{zetaprop}), $|\nabla \zeta| \equiv 0$ on $B(O,r_0/4)$, so
\begin{align*}
\int_{B} u_k^2 |\nabla \zeta|^2 = \int_{B\setminus B(O,r_0/4)} u_k^2
|\nabla \zeta|^2 \leq C,
\end{align*}
because $u \in C_{loc}^{1,1}(\Omega)$ and is therefore locally
bounded. Also, by (\ref{dul2}),
\begin{align*}
\int_{B} \zeta^2 |\nabla u_k|^2 \leq C
\end{align*}
independent of $k$.  From (\ref{pic1}) we conclude
\begin{align*}
\int_{B(O,r_0/4)} u_k^2 \leq C,
\end{align*}
and from the Monotone Convergence Theorem it follows
\begin{align*}
\int_{B} u_k^2 \leq C.
\end{align*}
\end{proof}

\begin{lemma} \label{weak} $u$ satisfies
\begin{align} \label{subsoln}
\Delta u \geq -\sigma_1(A)
\end{align}
in the $W^{1,2}$-sense on $B = B(O,r_0)$; i.e., for each
non-negative $\varphi \in C^1_{0}(B)$,
\begin{align} \label{weaksce}
\int - \langle \nabla u , \nabla \varphi \rangle  \geq \int -
\sigma_1(A) \varphi .
\end{align}
\end{lemma}

\begin{proof}
For $\epsilon > 0$ small, let $\eta_{\epsilon}$ denote the cut-off
function defined in (\ref{etadef}), and let $\varphi \in C_{0}^1(B)$
be non-negative.  Since (\ref{sce}) holds on $\Omega$ and
$\eta_{\epsilon}^2 \varphi$ is supported in $\Omega$ we have
\begin{align} \label{weak1'}
\int \eta_{\epsilon}^2 \varphi \Delta u \geq \int \ds
\frac{(n-2)}{2}|\nabla u|^2 \eta_{\epsilon}^2 \varphi - \int
\eta_{\epsilon}^2 \varphi \sigma_1(A).
\end{align}
Integrating by parts,
\begin{align*}
\int -\langle \nabla \varphi, \nabla u \rangle \eta_{\epsilon}^2 -
\int 2 \langle \nabla \eta_{\epsilon}, \nabla u \rangle
\eta_{\epsilon} \varphi \geq \int \ds \frac{(n-2)}{2}|\nabla u|^2
\eta_{\epsilon}^2 \varphi - \int \eta_{\epsilon}^2 \varphi
\sigma_1(A),
\end{align*}
which we rewrite as
\begin{align} \label{weak2}
\int -\langle \nabla \varphi, \nabla u \rangle \eta_{\epsilon}^2
\geq \int \Big[ 2 \langle \nabla \eta_{\epsilon}, \nabla u \rangle
\eta_{\epsilon} \varphi + \ds \frac{(n-2)}{2}|\nabla u|^2
\eta_{\epsilon}^2 \varphi \Big] - \int \eta_{\epsilon}^2 \varphi
\sigma_1(A).
\end{align}
By the Lebesgue Dominated convergence theorem, for the first and
last integrals in (\ref{weak2}) we have
\begin{align} \label{lims} \begin{split}
\int -\langle \nabla \varphi, \nabla u \rangle \eta_{\epsilon}^2
&\to \int -\langle \nabla \varphi, \nabla u \rangle, \quad \mbox{and} \\
\\
\int \eta_{\epsilon}^2 \varphi \sigma_1(A) &\to \int \varphi
\sigma_1(A) \end{split}
\end{align}
as $\epsilon \to 0$.  We estimate the middle integral in the
following way:  First,
\begin{align*}
\int 2 \langle \nabla \eta_{\epsilon}, \nabla u \rangle
\eta_{\epsilon} \varphi \geq \int - \ds \frac{(n-2)}{2}|\nabla u|^2
\eta_{\epsilon}^2 \varphi - \int \ds \frac{2}{(n-2)}\varphi |\nabla
\eta_{\epsilon}|^2.
\end{align*}
Therefore,
\begin{align} \label{weak3}
\int \Big[ 2 \langle \nabla \eta_{\epsilon}, \nabla u \rangle
\eta_{\epsilon} \varphi + \ds \frac{(n-2)}{2}|\nabla u|^2
\eta_{\epsilon}^2 \varphi \Big] \geq \int - C_n \varphi |\nabla
\eta_{\epsilon}|^2.
\end{align}
Note that
\begin{align} \label{lims2}
\int  \varphi |\nabla \eta_{\epsilon}|^2 &\leq C \epsilon^{-2}
\int_{B(O,2\epsilon) \setminus B(O,\epsilon)} \varphi \\
&\leq C\epsilon^{(n-2)} \to 0
\end{align}
as $\epsilon \to 0$.  Substituting (\ref{lims}),(\ref{weak3}), and
(\ref{lims2}) into (\ref{weak2}), we get (\ref{weaksce}).
\end{proof}

To complete the proof of the Proposition, we refer to Theorem
8.17 of Gilbarg-Trudinger \cite{GT}, which implies that any
$W^{1,2}$-solution of (\ref{subsoln}) satisfies
\begin{align*}
\sup_{B(O,r_0/2)} u \leq C(r_0,g) \Big( \| u \|_{L^2(B(O,r_0))} +
C(g) \Big).
\end{align*}
Thus, the desired bound follows from Lemma \ref{uinw12}.  $\Box$

\end{proof}

While Proposition \ref{prop1} gives an upper bound on solutions,
there are examples of singular solutions for which a lower bound
fails to hold. The next result controls the rate at which $u(x)$ can
go to $-\infty$, provided $u$ satisfies (\ref{gradr1}).

\begin{proposition} \label{prop2}
Let $u \in C^{1,1}_{loc}(\Omega)$ satisfy
\begin{align*}
 \sigma_1(A + \nabla^2 u
+ du \otimes du - \frac{1}{2}| \nabla u|^2 g) \geq 0
\end{align*}
almost everywhere in $\Omega = B(O,r_0) \setminus \{O\}$.  Assume
$u$ satisfies the scale-invariant gradient estimate (\ref{gradr1}).
Then there is a constant $C>0$ such that
\begin{align} \label{lwbd}
u(x) \geq 2 \log d(x) - C.
\end{align}
\end{proposition}


The Proof of Proposition \ref{prop2} is essentially contained in
Proposition 6.1 of \cite{GVRic}.  The only difference is the
regularity assumed: we need to show the same argument applies to
$C^{1,1}$-solutions.

\begin{proof}
As we observed above, $u \in C_{loc}^{1,1}(\Omega)$ satisfies the
inequality
\begin{align*}
\Delta u \geq \ds \frac{(n-2)}{2}|\nabla u|^2 - \sigma_1(A)
\end{align*}
a.e. in $\Omega$.  Let
\begin{align} \label{wdef}
w =  e^{- \frac{(n-2)}{2}u} .
\end{align}
Then $w \in C_{loc}^{1,1}(\Omega)$, and a simple calculation shows
that $w$ satisfies
\begin{align} \label{sce2}
Lw = \Delta w - \frac{(n-2)}{4(n-1)}Rw \leq 0
\end{align}
a.e. in $\Omega$.  Let $\Gamma$ denote the Green's function for $L$
with pole at $O$.  Since $\Gamma(x) \sim d(x)^{2-n}$, it suffices to
prove
\begin{align} \label{goal1}
w(x) \leq C\Gamma(x), \quad x \in \Omega
\end{align}
for some constant $C$.  To this end, consider the function
\begin{align} \label{Gdef}
G(x) = \ds \frac{w(x)}{\Gamma(x)}.
\end{align}
It follows from (\ref{sce2}) and the definition of $\Gamma$ that
\begin{align} \label{GDE}
\Delta G \leq - 2 \langle \nabla G, \ds \frac{\nabla \Gamma}{\Gamma}
\rangle
\end{align}
a.e. in $\Omega$.  Since clearly $G \in C_{loc}^{1,1}(\Omega)$, it
follows that (\ref{GDE}) holds in a $W^{1,2}$-sense in $\Omega$.

Now, fix $r>0$ small, and let $\Omega_r = B \setminus B(O,r)$. By
the strong maximum principle (\cite{GT}, Theorem 8.19), $G$ cannot
attain an interior minimum in $\Omega_r$ unless it is constant; of
course, if $G$ were constant then (\ref{goal1}) would follow
immediately.  Therefore, assume $G$ attains its minimum on $\partial
\Omega_r = \partial B(O,r_0) \cup \partial B(O,r)$; in fact, assume
it is attained on $\partial B(O,r)$.

Since we are assuming $u$ satisfies a scale-invariant gradient
estimate, given any two points $x,y \in \partial B(O,r)$ we have
\begin{align*}
|u(x) - u(y)| \leq Cr \max_{\partial B(O,r)} |\nabla u| \leq C.
\end{align*}
From this inequality it follows that $w$--and hence $G$--satisfies a
Harnack inequality:
\begin{align} \label{har1}
\max_{\partial B(O,r)} G \leq C \min_{\partial B(O,r)} G.
\end{align}
Therefore,
\begin{align*}
\min_{\Omega_r} G = \min_{\partial \Omega_r} G \geq C^{-1}
\max_{\partial B(O,r)} G.
\end{align*}
In case the minimum of $G$ is attained on $\partial B(O,r_0)$, we
can apply the same argument.  In either case, we conclude
\begin{align} \label{har2}
\max_{\partial \Omega_r} G \leq C \min_{\Omega_r} G.
\end{align}
If we choose a point $x_0 \in \Omega$, then (\ref{har2}) implies
\begin{align}  \label{har3}
\max_{\partial \Omega_r} G \leq C \min_{\Omega_r} G \leq C G(x_0) =
C^{\prime},
\end{align}
Since (\ref{har3}) holds for all $r>0$ small, it follows that $G$ is
uniformly bounded in $\Omega$.  This completes the proof.
\end{proof}

\subsection{An integral estimate}

The final result of this section is an integral estimate which will
be used in the proof of Theorem \ref{thm1'}.

\begin{proposition} \label{Lnest}
Let $u \in C^{1,1}_{loc}$ satisfy $A_u \in
\overline{\Gamma}_{\delta}$ almost everywhere in $\Omega = B(O,r_0)
\setminus \{ O \}.$ Let
\begin{align*}
v = \ds e^{\beta u},
\end{align*}
where $\beta$ is defined in (\ref{betadef}).
Then $v$ satisfies
\begin{align} \label{LnCor}
\int_{\Omega} |\nabla v|^n \leq C(\delta,n,g).
\end{align}
\end{proposition}

\begin{proof} This estimate is actually a corollary
of Theorem 3.7 of \cite{GVRic}:

\begin{theorem}(\cite[Theorem 3.7]{GVRic})
\label{GVRicest}
Let $u \in
C_{loc}^{1,1}\big(A(\frac{1}{2}r_1,2r_2)\big)$, where $O \in M^n$
and $A(\frac{1}{2}r_1,2r_2)$ denotes the annulus
$A(\frac{1}{2}r_1,2r_2) \equiv B(O,2r_2) \setminus \overline{B(O,
\frac{1}{2}r_1)}$, with $0 < r_1 < r_2$. Assume $g_u = e^{-2u}g$
satisfies
\begin{align} \label{posRicci}
Ric(g_u) - 2\delta_0 \sigma_1(A_u)g \geq 0
\end{align}
almost everywhere in $A(\frac{1}{2}r_1,2r_2)$ for some $0 \leq
\delta_0 < \frac{1}{2}$. Define
\begin{align} \label{alphadeltadef}
\alpha_{0} = \frac{(n-2)}{(1-2\delta)}\delta_0 \geq 0.
\end{align}
Then given any $\alpha > \alpha_{0}$, there are constants $p \geq n$
and $C = C((\alpha - \alpha_{0})^{-1},n) > 0$ such that
\begin{align} \label{localintgrad0} \begin{split}
\int_{A(r_1,r_2)} |\nabla u|^p e^{\alpha u} dvol_g &  \leq C\Bigg(
\int_{A(\frac{1}{2}r_1,2r_2)} |Ric_g|^{p/2} e^{\alpha u} dvol_g \\
&+ r_1^{-p} \int_{A(\frac{1}{2}r_1,r_1)}e^{\alpha u}dvol_g +
r_2^{-p} \int_{A(r_2, 2r_2)} e^{\alpha u} dvol_g \Bigg).
\end{split}
\end{align}
In fact, we can take
\begin{align} \label{pdef}
p = n + 2\alpha_{0} \geq n.
\end{align}
\end{theorem}

\begin{remark}  Keeping with the conventions of this paper we will
omit the volume form $dvol_g$.
\end{remark}
Let $\delta_0$ and $\delta$ be related by the formula
\begin{align}
\delta = \frac{ 1 - 2 \delta_0}{n-2},
\end{align}
or equivalently,
\begin{align}
\label{delt0def}
2 \delta_0 = 1 + (2-n) \delta.
\end{align}
Inequality (\ref{posRicci}) for $\delta_0$ is equivalent to saying
that  $A_u \in \overline{\Gamma}_{\delta}$.
Since $A_u \in \overline{\Gamma}_{\delta}$ almost everywhere in
$\Omega = B(O,r_0) \setminus \{ O \}$, inequality (\ref{posRicci})
holds for $\delta_0$ as defined in (\ref{delt0def}).
Clearly, the inequality (\ref{posRicci}) then also holds for
$\delta_0 = 0$, and hence $\alpha_0 = 0$ and
$p = n$. Letting $r_2 = \frac{1}{2}r_0$ and $r_1 = r <
\frac{1}{2}r_0$, from (\ref{localintgrad0}) we have
\begin{align} \label{localintgrad1} \begin{split}
\int_{A(r,r_0/2)} |\nabla u|^n e^{\alpha u} &  \leq
C(\alpha^{-1},r_0,g) \Bigg(
\int_{A(r/2,r_0)} e^{\alpha u} \\
&+ r^{-n} \int_{A(r/2,r)}e^{\alpha u} + \int_{A(r_0/2, r_0)}
e^{\alpha u}  \Bigg),
\end{split}
\end{align}
for any $\alpha > 0$.  By Proposition \ref{prop1}, $u$ is bounded
above on $\Omega$, and therefore $e^{\alpha u} \leq C$ on $\Omega$.
Also, notice the middle integral on the right-hand side of
(\ref{localintgrad1}) is uniformly bounded:
\begin{align*}
r^{-n} \int_{A(r/2,r)}e^{\alpha u} &\leq Cr^{-n}
\int_{A(r/2,r)}  \\
&\leq Cr^{-n} \big( cr^n \big) \\
&\leq C.
\end{align*}
Consequently, for all $\alpha > 0$ we have
\begin{align} \label{localintgrad2goog}
\int_{A(r,r_0/2)} |\nabla u|^n e^{\alpha u} & \leq
C(\alpha^{-1},r_0,g)
\end{align}
independent of $r$.  Letting $r \rightarrow 0$ we obtain
\begin{align} \label{localintgrad2}
\int_{\Omega} |\nabla u|^n e^{\alpha u} &  \leq
C(\alpha^{-1},r_0,g).
\end{align}
If $v = e^{\beta u}$, then
\begin{align*}
|\nabla v|^n = \beta^n |\nabla u|^n e^{\beta n u}.
\end{align*}
Therefore, taking $\alpha = \beta n$ in (\ref{localintgrad2}) we get
(\ref{LnCor}).

\end{proof}

\section{The proof of Theorems \ref{thm1} and \ref{thm1'}}

To prove Theorems \ref{thm1} and \ref{thm1'} we assume $u \in
C^{1,1}_{loc}(\Omega)$ with $A_u \in \overline{\Gamma}_{\delta}$ almost
everywhere. For part $(i)$, we further assume
\begin{align} \label{goodinf2}
\liminf_{x \to O} u >  -\infty,
\end{align}
and that either $g$ is $LCF$, or that $u$ satisfies a scale-invariant
$C^2$-estimate. In each case we wish to show that $u$ can be
extended to a H\"older continuous function $u^{\ast} \in
C^{1,1}_{loc}(\Omega) \cap C^{\gamma}(B(O,r_0))$, where
$\gamma$ is defined in (\ref{gammacon}).
Note that (\ref{goodinf2}) along with (\ref{supbound}) imply
\begin{align} \label{linftyu}
\sup_{\Omega} |u| < \infty.
\end{align}
As in the proof of Theorem \ref{adtocon}, let $v = e^{\beta u}$,
where $\beta$ is defined in (\ref{betadef}).
By (\ref{linftyu}), $v \in C^{1,1}_{loc}(\Omega)$ satisfies
\begin{align} \label{L0v}
0 < c_0 \leq v(x) \leq c_0^{-1}.
\end{align}

First, consider the locally conformally flat case.   Let $\{ x^i
\}$ denote conformally flat coordinates centered at $O$; in this
coordinate system (\ref{hessvcone}) is equivalent to
\begin{align} \label{localcone2}
\left( \partial_i \partial_j v + \beta v A_{ij} \right) \in
\overline{\Gamma}_{\delta}
\end{align}
a.e. in $\Omega$.  Define
\begin{align}
W(x) = v(x) + \Lambda|x|^2,
\end{align}
where $|x|^2 = \sum_i (x^i)^2$ and $\Lambda >> 0$ is a large
constant.  Then
\begin{align*}
\partial_i \partial_j W &= \partial_i \partial_j v + 2\Lambda
\delta_{ij} \\
&= \partial_i \partial_j v + \beta v A_{ij} + 2\Lambda \delta_{ij} - \beta v A_{ij}.
\end{align*}
Therefore, by (\ref{L0v}) and (\ref{localcone2}), we can choose $\Lambda >> 0$ large enough so that
\begin{align}
\label{Vcon2}
\partial_i \partial_j W \in \Gamma_{\delta}
\end{align}
a.e., in a deleted neighborhood of $O$.

Using the coordinates $\{ x^i \}$ we can identify a
neighborhood of $O \in M^n$ with a neighborhood $U$ of the origin $0
\in \mathbf{R}^n$.  Following \cite{TrudingerWangII}, let $\rho \in
C_0^{\infty}$ be a spherically symmetric mollifier satisfying
$\rho(x) \geq 0$, $\rho(x) = 1$ for $|x| < 1$, $\rho(x) = 0$ for
$|x|
> 2$, and $\int \rho = 1$.  Define the mollification of $W$ by
\begin{align} \label{Vhdef}
W_h(x) = h^{-n} \int \rho(\frac{x-y}{h}) W(y) dy.
\end{align}
Let $U'$ be a subset of $U$, such that the $h$-neighborhood of
$U'$ is also contained in $U$.

\begin{proposition}
$W_h : U' \rightarrow \mathbf{R}$ is a smooth, bounded, strictly
$\delta$-convex function.
\end{proposition}

\begin{proof}
The smoothness of $W$ follows from elementary properties of
convolutions. We let $x \in U$ with $d(x, \partial U) > h$ and $r
>0$ a small number.  By the divergence theorem,
\begin{align*}
D_{ij} W_h(x) &= \int_{U} D_{ij}\rho_h(x-y) W(y) dy\\
& = \int_{B(0,r)} D_{ij}\rho_h(x-y) W(y) dy
+  \int_{U \setminus B(0,r)}  D_{ij}\rho_h(x-y) W(y) dy\\
& =  \int_{B(0,r)} D_{ij}\rho_h(x-y) W(y) dy
- \int_{U \setminus B(0,r)}  D_{i}\rho_h(x-y) D_j W(y) dy\\
& \ \ \ +  \oint_{ \partial B(0,r)} \nu_j D_{i}\rho_h(x-y) W(y)
dS(y),
\end{align*}
where $\{ \nu_j \}$ are the components of the outward unit normal
to $\partial B(0,r).$ Since $W$ is bounded, as $r \rightarrow 0$
we obtain
\begin{align}
D_{ij} W_h(x) &= - \int_{U}  D_{i}\rho_h(x-y) D_j W(y) dy.
\end{align}
Applying the divergence theorem again,
\begin{align} \label{star} \begin{split}
D_{ij} W_h(x) &= -\int_{B(0,r) }  D_{i}\rho_h(x-y) D_j W(y) dy
- \int_{U \setminus B(0,r) }  D_{i}\rho_h(x-y) D_j W(y) dy\\
& = - \int_{B(0,r) }  D_{i}\rho_h(x-y) D_j W(y) dy
+ \int_{U \setminus B(0,r) } \rho_h(x-y) D_{ij} W(y) dy\\
& \ \ \ \ \ - \oint_{ \partial B(0,r)}   \rho_h(x-y)  D_j W(y)
dS(y). \end{split}
\end{align}

\begin{lemma}
\label{lcflemma}
There is a sequence $r_i \rightarrow 0$ such that
\begin{align} \label{intr1}
\Big| \oint_{ \partial B(0,r_i)}  \rho_h(x-y) \nu_i D_j W(y) dS(y)
\Big|  \rightarrow 0,
\end{align}
\begin{align} \label{intr2}
\Big| \int_{B(0,r_i) }  D_{i}\rho_h(x-y) D_j W(y) dy \Big| \rightarrow 0,
\end{align}
as $i \to \infty$.
\end{lemma}

\begin{proof}

By Proposition \ref{Lnest}, $v$ satisfies
\begin{align} \label{LnCor1}
\int_{\Omega} |\nabla v|^n \leq C(\delta,n,g).
\end{align}
This implies
\begin{align} \label{LnW}
\int_{U} |DW(y)|^n dy \leq C.
\end{align}
It follows from the co-area formula that there is a sequence of radii $r_i \rightarrow 0$ such that
\begin{align} \label{slice}
\oint_{\partial B(0,r_i)} | DW(y)|^n dS(y) \leq C/ r_i.
\end{align}
Therefore, by H\"older's inequality,
\begin{align*}
\Big| \oint_{ \partial B(0,r_i)}  \rho_h(x-y) \nu_i D_j W(y) dS(y) \Big| &\leq C \oint_{ \partial B(0,r_i)} |DW(y)| dS(y) \\
& \leq C r_i^{(n-1)^2/n} \Big( \oint_{ \partial B(0,r_i)} |DW(y)|^n dS(y) \Big)^{1/n} \\
& \leq C r_i^{n-2},
\end{align*}
and (\ref{intr1}) follows.

For the same sequence of radii, using (\ref{LnW}) and H\"older's inequality we have
\begin{align*}
\Big| \int_{B(0,r_i) }  D_{i}\rho_h(x-y) D_j W(y) dy \Big| & \leq C \int_{B(0,r_i) } |DW(y)| dy \\
&\leq Cr_i^{n-1}\Big( \int_{B(0,r_i) } |DW(y)|^n dy \Big)^{1/n} \\
& \leq Cr_i^{n-1}.
\end{align*}
\end{proof}

Taking $r = r_i$ in (\ref{star}) and letting $i \rightarrow \infty$ we obtain
\begin{align}
\label{qeust} D_{ij} W_h(x) = \int_{U} \rho_h(x-y) D_{ij} W(y) dy.
\end{align}
Since $D^2 W \in \Gamma_{\delta}$ almost everywhere in $U$, it follows that $D_{ij}W_h(x) \in
\Gamma_{\delta}$.  Consequently, $W_h$ is a smooth, strictly
$\delta$-convex function on $U'$.
\end{proof}

\begin{proposition} $W$ has a $C^{\gamma}$-H\"older continuous
extension across the origin.
\end{proposition}
\begin{proof}
From Theorem \ref{unifhold} and the Arzela-Ascoli Theorem,
$W_{h_i} \rightarrow \overline{W}$ uniformly in the $C^{\gamma'}$-norm,
$\gamma' < \gamma$, for some sequence $h_i \rightarrow 0$, and
$\overline{W} \in C^{\gamma}(B)$.  Define $W(0) =
\overline{W}(0)$. From general properties of mollification, $W_{h_i}
\rightarrow W$ on $B \setminus \{0 \}$. Therefore $W = \overline{W}$
in $B$.
\end{proof}

Since $v = W - \Lambda|x|^2$, the same holds for $v$ and (by the definition of $v$) for $u$ as well.
\vskip.1in

Turning to the non-$LCF$ case, we now assume $u$ satisfies a scale-invariant $C^2$-estimate.
In view of (\ref{linftyu}) and inequality (\ref{gradr}), $v$ satisfies
\begin{align} \label{C1v}
|\nabla^2v(x)| + |\nabla v(x)|^2 \leq c_2d(x)^{-2}
\end{align}
for almost every $x \in \Omega$.

Let $\{ x^i \}$ denote normal coordinates (with
respect to the background metric $g$) centered at $O$.  In this
coordinate system (\ref{hessvcone}) is equivalent to
\begin{align} \label{localcone3}
g^{jl}\left( \partial_i \partial_j v - \Gamma_{ij}^k \partial_k v +
\beta v A_{ij} \right) \in \overline{\Gamma}_{\delta}
\end{align}
a.e. in $\Omega$. In normal coordinates, $g^{jl} = \delta^{jl} +
O(|x|^{2}),$ and $|\Gamma_{ij}^k| = O(|x|)$ as $|x| \rightarrow 0$,
so using (\ref{L0v}) and (\ref{C1v}) we conclude
\begin{align}
\label{localcone4}
\partial_i \partial_j v + B_{ij} \in \overline{\Gamma}_{\delta},
\end{align}
where $B_{ij}$ is uniformly bounded:
\begin{align} \label{Bbound}
\Vert B_{ij} \Vert < C_1.
\end{align}
As before, we let
\begin{align*}
W(x) = v(x) + \Lambda |x|^2.
\end{align*}
Then
\begin{align*}
\partial_i \partial_j W &= \partial_i \partial_j v + 2\Lambda
\delta_{ij} \\
&= \partial_i \partial_j v + B_{ij} + 2\Lambda \delta_{ij} - B_{ij}.
\end{align*}
Therefore, by (\ref{localcone4}) and (\ref{Bbound}) we can choose $\Lambda >> 0$ large enough so that
\begin{align}
\label{Vcon2goog}
\partial_i \partial_j W \in \Gamma_{\delta}
\end{align}
a.e., in a deleted neighborhood of $O$.
The rest of the proof proceeds exactly as in the $LCF$-case (note
however that Lemma \ref{lcflemma} is much easier to prove
under the assumption (\ref{gradr1})).
This completes the proof of part $(i)$ of Theorem \ref{thm1} and
Theorem \ref{thm1'}.


\vskip.2in

To prove part $(ii)$ of Theorems \ref{thm1} and \ref{thm1'}, we assume
\begin{align} \label{goodinf2b}
\liminf_{x \to O} u =  -\infty,
\end{align}
and first assume that either $g$ is $LCF$; or that
$u$ satisfies (\ref{gradr}). In each case we wish to show that $u$ obeys the growth estimate
\begin{align} \label{sharpie}
u(x) \leq 2 \log d(x) + C.
\end{align}

To this end, we appeal to part $(i)$, in which we showed that the
function $v = e^{\beta u}$ can be extended to a H\"older continuous
function $v^{\ast} \in C^{1,1}_{loc}(\Omega) \cap
C^{\gamma}(B(O,r_0))$, where $\gamma$ is given by (\ref{gammacon}).  Therefore,
\begin{align*}
|v^{\ast}(x) - v^{\ast}(y)| \leq C d(x,y)^{\gamma}
\end{align*}
for all $x,y$ near $O$.
Rewriting this in terms of $u$, we have
\begin{align*}
| e^{\beta u(x)} - e^{ \beta u(y)}| \leq C d(x,y)^{\gamma}
\end{align*}
for all $x,y \in \Omega$.
Since $\beta = \gamma/2$, this implies
\begin{align} \label{holde}
| e^{u(x)} - e^{u(y)}| \leq C d(x,y)^{2}.
\end{align}
From assumption (\ref{goodinf2b}), there exists a sequence of
points $y_i \in \Omega$ with
$y_i \rightarrow O$ and $u(y_i) \rightarrow -\infty$ as $i
\rightarrow \infty$.  Taking $y = y_i$ in (\ref{holde}) and letting
$i \rightarrow \infty$ we obtain
\begin{align*}
e^{u(x)} \leq C d(x)^{2},
\end{align*}
which implies
\begin{align} \label{above}
u(x) \leq 2 \log d(x) + C.
\end{align}
To complete the proof of Theorems \ref{thm1} and \ref{thm1'}, we
turn our attention to the lower inequality (\ref{Greensgrowth}).
Since in both cases we are assuming $u$ satisfies scale-invariant
$C^1$-estimate it follows from Proposition \ref{prop2} that
\begin{align} \label{oneside}
u(x) \geq 2 \log d(x) - C.
\end{align}

\begin{remark} The method in this section simplifies somewhat our
previous proof of the above estimate (\ref{above}) which was given in Section $6$ of
\cite{GVRic}. However, the methods are in essence the same in that
they are both based on some version of the maximum principle.
\end{remark}

\section{The scale-invariant estimates}
\label{Pointwise}

In this section we verify the scale-invariant estimates
(\ref{gradr1}) for $C^3$-solutions and (\ref{gradr}) for
$C^4$-solutions subject to various assumptions.  These results
follow from various local estimates for solutions, along with a
scaling argument.

The local gradient and $C^2$-estimates
for (\ref{hesskF}) and (\ref{rats}) were first proved in
\cite{GuanWang1} and \cite{GuanLinWang2}. This work was
extended to local gradient and $C^2$-estimates for
(\ref{sigmakt}) in \cite{LiLi2}.
For the case of general symmetric functions $F$ satisfying
assumptions ${ \bf (i)-(iv)}$ in the introduction, local $C^2$-estimates
for $C^4$ solutions were
proved recently in \cite{Sophie1}. We note that the local gradient
estimates for $C^3$ solutions are currently not known for general $F$,
but we conjecture these to be true.

We say that $u \in C^2$ is $k-$admissible if $A_u \in
\overline{\Gamma_{\sigma_k}}$ (this is not to be confused
with the notion of $\delta$-admissibility defined above).
For equations (\ref{hesskF}) and
(\ref{rats}), the scale-invariant estimates are based on the following
results of Guan-Wang and Guan-Lin-Wang:

\begin{theorem} {\em (Theorem 1.1 of \cite{GuanWang1}, Theorem 1 of \cite{GuanLinWang2})} \label{GWlocal}
Let $u \in C^3(M^n)$ be a $k-$admissible solution of (\ref{hesskF}) or (\ref{rats})
in $B(x_0,\rho)$, where $x_0 \in M^n$ and $\rho
> 0$.  Then there is a constant
\begin{align*}
C_0 = C_0(k,n,\rho,\|g\|_{C^2(B(x_0,\rho))},
\|f\|_{C^1(B(x_0,\rho))}),
\end{align*}
such that
\begin{align} \label{GWlocalest}
|\nabla u|^2(x) \leq C_0\big(1 + e^{-2 \inf_{B(x_0,\rho)}u}\big)
\end{align}
for all $x \in B(x_0,\rho/2)$.

Let $u \in C^4(M^n)$ be a $k-$admissible solution of (\ref{hesskF})
in $B(x_0,\rho)$, where $x_0 \in M^n$ and $\rho
> 0$.  Then there is a constant
\begin{align*}
C_0 = C_0(k,n,\rho,\|g\|_{C^3(B(x_0,\rho))},
\|f\|_{C^2(B(x_0,\rho))}),
\end{align*}
such that
\begin{align} \label{GWlocalest2}
|\nabla^2 u|(x) + |\nabla u|^2(x) \leq C_0\big(1 + e^{-2
\inf_{B(x_0,\rho)}u}\big)
\end{align}
for all $x \in B(x_0,\rho/2)$.

\end{theorem}

We first observe that when $f(x) \equiv 0$ in (\ref{hesskF}) or (\ref{rats}),
then there is no exponential term in the estimate (\ref{GWlocalest}).
We will only verify this explicitly for solutions of (\ref{hesskF}),
but the argument for solutions of (\ref{rats}) is essentially identical.

\begin{corollary} \label{GWCorf0} Let $u \in C^3(M^n)$ be a $k-$admissible solution of (\ref{hesskF})
in $B(x_0,\rho)$ with $f(x) \equiv 0$. Then there is a constant
\begin{align*}
C_0 = C_0(k,n,\rho,\|g\|_{C^2(B(x_0,\rho))})
\end{align*}
such that
\begin{align} \label{GWlocalestf0}
|\nabla u|^2(x) \leq C_0
\end{align}
for all $x \in B(x_0,\rho/2)$. In fact,
\begin{align} \label{Crad}
C_0 = C_1 \rho^{-2},
\end{align}
where $C_1 = C_1(k,n,\|g\|_{C^2(B(x_0,\rho))})$.

 Let $u \in C^4(M^n)$ be a $k-$admissible solution of (\ref{hesskF})
in $B(x_0,\rho)$ with $f(x) \equiv 0$. Then there is a constant
\begin{align*}
C_0 = C_0(k,n,\rho,\|g\|_{C^3(B(x_0,\rho))})
\end{align*}
such that
\begin{align} \label{GWlocalestf1}
|\nabla^2 u|(x) + |\nabla u|^2(x) \leq C_0
\end{align}
for all $x \in B(x_0,\rho/2)$. In fact,
\begin{align} \label{Crad2}
C_0 = C_1 \rho^{-2},
\end{align}
where $C_1 = C_1(k,n,\|g\|_{C^3(B(x_0,\rho))})$.

\end{corollary}

\begin{proof} If we
imitate the proof of Guan and Wang, one can trace the origin of the
exponential term in (\ref{GWlocalest}) to two places: inequalities
(2.20) and (3.10) in \cite{GuanWang1}.  These inequalities appear
when estimating the gradient and Hessian terms respectively.

For the gradient term, Guan and Wang estimate
\begin{align*}
T_1 = \sum_{l} F_l u_l &= \sum_{l}
\big(fe^{-2u}\big)_l u_l \\
&= \sum_{l} e^{-2u}(f_lu_l - 2fu_l^2),
\end{align*}
where the subscript $l$ denotes $\frac{\partial}{\partial x_l}$.
They do so in the following way:
\begin{align*}
T_1 = \sum_{l} F_l u_l &= \sum_{l} e^{-2u}(f_lu_l - 2fu_l^2)\\
&\geq -C(1+e^{-2u})|\nabla u|^2,
\end{align*}
thus the appearance of the exponential term in (\ref{GWlocalest}).
Of course, if $f \equiv 0$ then $T_1 \equiv 0$.  While estimating
the Hessian a similar term appears in (3.10) of \cite{GuanWang1}:
\begin{align*}
T_2 = \sum_{l} \rho F_{ll},
\end{align*}
where $\rho$ is a cut-off function supported in $B(0,r_0)$.
Unfortunately, a typographical appears when estimating $T_2$; the
first line of (3.10) should read
\begin{align*}
T_2 = \sum_{l} \rho F_{ll} &= \sum_{l} \rho \big( fe^{-2u}
\big)_{ll} \\
&= \sum_{l} \rho \big( f_{ll} - 4 f_{l}u_{l} + 4 f u^2_{l} - 2f
u_{ll} \big) e^{-2u}.
\end{align*}
In any case, once again if $f \equiv 0$ then $T_2 \equiv 0$, and the
Hessian estimate no longer contributes an exponential term.
Consequently, inequality (\ref{GWlocalestf1}) holds.

To prove (\ref{Crad2}) we need to specify the dependence of $C_0$ on
the radius of the ball $\rho$; that is, we need to show
\begin{align}  \label{sigef}
|\nabla^2 u|(x) + |\nabla u|^2(x) \leq C_1(g)/\rho^2.
\end{align}
However, such an estimate follows from an elementary scaling
argument; see Lemma 3.3 of \cite{GVRic} for a detailed explanation.

\end{proof}

From Corollary \ref{GWCorf0} it follows that the scale-invariant
estimate (\ref{gradr}) holds for solutions of (\ref{hesskF}) when $f
\equiv 0$:

\begin{corollary} \label{corf0}
Let $u \in C^3(\Omega)$ be a (degenerate) admissible solution of
(\ref{hesskF}) in $\Omega = B(O,r_0) \setminus \{O\}$ with $f \equiv
0$ in a neighborhood of $O$.  Then $u$ satisfies (\ref{gradr1}).

Let $u \in C^4(\Omega)$ be a (degenerate) admissible solution of
(\ref{hesskF}) in $\Omega = B(O,r_0) \setminus \{O\}$ with $f \equiv
0$ in a neighborhood of $O$.  Then $u$ satisfies (\ref{gradr}).
\end{corollary}

\begin{proof}  For $x \in \Omega$ close enough to $O$, $u$ is a
solution of (\ref{hesskF})
in the ball $B = B(x, d(x)/2)$, where $d(x) = d(O,x)$. By Corollary
\ref{GWCorf0}, if $u$ is $C^3$, then $u$ satisfies
\begin{align*}
|\nabla u|^2(x) \leq Cd(x)^{-2},
\end{align*}
and if $u$ is $C^4$, then $u$ satisfies
\begin{align*}
|\nabla^2 u|(x) + |\nabla u|^2(x) \leq Cd(x)^{-2}.
\end{align*}
\end{proof}

\subsection{The Proof of Proposition \ref{prop3}} 
If one traces through the local estimates of Guan-Lin-Wang for the
quotient equations (\ref{rats}) (\cite{GuanLinWang2}), and in
\cite{GuanWangDuke}, or the estimates for solutions of
(\ref{sigmakt}) in \cite{LiLi2}, or the estimates for general
$F$ in \cite{Sophie1}, then in all cases the presence of the
exponential term comes from the right-hand side of the equation.
Therefore, when $f \equiv 0$, the same argument presented above
leads to the scale-invariant estimates for solutions of
these equations.  $\Box$

\section{H\"older Extension}

Next we consider the case of finite volume metrics.  As a
preliminary observation, we note that a corollary of the local
estimates is an $\epsilon$-regularity result:

\begin{proposition} \label{epsilonreg} {\em (See Proposition 3.6
of \cite{GuanWang1} and Proposition 3.4 of \cite{GVRic})}
\vskip.1in
\noindent $(i)$ Let $u \in C^3(B(x_0,\rho))$ be an admissible
solution of (\ref{hesskF}), (\ref{rats}), or (\ref{sigmakt}). Then
there exist constants $\epsilon_0
> 0$ and $C = C(\epsilon_0)$ such that if
\begin{align} \label{smallenergy}
\int_{B(x_0,\rho)} e^{-nu} dvol_g \leq \epsilon_0,
\end{align}
then
\begin{align} \label{lowerb1}
\inf_{B(x_0,\rho/2)} u \geq -C_2 + \log \rho.
\end{align}  Consequently, there is a constant
\begin{align*}
C_3 = C_3(k,n,\epsilon_0,\|g\|_{C^3(B(x_0,\rho))}),
\end{align*}
such that
\begin{align} \label{gradonr1}
|\nabla u|^2(x) \leq C_3\rho^{-2}
\end{align}
for all $x \in B(x_0,\rho/4)$.
\vskip.1in
\noindent $(ii)$ Let $u \in C^4(B(x_0,\rho))$ be an admissible
solution of (\ref{hessF}). Then
there exist constants $\epsilon_0
> 0$ and $C = C(\epsilon_0)$ such that if
\begin{align} \label{smallenergy2}
\int_{B(x_0,\rho)} e^{-nu} dvol_g \leq \epsilon_0,
\end{align}
then
\begin{align} \label{lowerb2}
\inf_{B(x_0,\rho/2)} u \geq -C_2 + \log \rho.
\end{align}  Consequently, there is a constant
\begin{align*}
C_4 = C_4(k,n,\epsilon_0,\|g\|_{C^4(B(x_0,\rho))}),
\end{align*}
such that
\begin{align} \label{gradonr2}
|\nabla^2 u|(x) +  |\nabla u|^2(x) \leq C_4\rho^{-2}
\end{align}
for all $x \in B(x_0,\rho/4)$.

\end{proposition}

\begin{proof}  For solutions of (\ref{hesskF}), (\ref{rats}), and
(\ref{sigmakt}), this was proved in Proposition 3.6
of \cite{GuanWang1} and Proposition 3.4 of \cite{GVRic}.
For $C^4$ solutions of (\ref{hessF}), this follows from
the same method of proof, using the local $C^2$-estimates
from \cite{Sophie1}.
\end{proof}
\subsection{The Proof of Theorems \ref{thm5} and \ref{thm5'}} 
Suppose
$u \in C^4$ satisfies the hypotheses of Theorem \ref{thm5}. Since
the volume of $g_u = e^{-2u}g$ is finite,
\begin{align}  \label{volup}
Vol(g_u) = \int_{\Omega} e^{-nu} < \infty.
\end{align}
Therefore, there is a radius $r_1 = r_1(\epsilon_0)$ such that for
all $x \in \Omega$ with $d(x) < r_1$,
\begin{align*}
\int_{B(x,d(x)/2)} e^{-nu} < \epsilon_0.
\end{align*}
From (\ref{lowerb2}) and (\ref{gradonr2}) it follows that
\begin{align} \label{slowgro}
u(x) \geq -C + \log d(x), \quad
 |\nabla^2 u(x)| + |\nabla u(x)|^2 \leq Cd(x)^{-2}.
\end{align}
This proves Theorem \ref{thm5}. In the case that $u \in C^3$,
and $F$ is as in cases (\ref{hesskF}), (\ref{rats}), or (\ref{sigmakt}),
the above argument yields the scale
invariant gradient estimate (\ref{gradr1}), which
proves Theorem \ref{thm5'}. \quad \quad $\Box$
\subsection{The Proof of Theorem \ref{thm6}}
\label{fgeq0}

In Theorem \ref{thm6} we assume that $g$ is $LCF$ and $u \in C^3$ is a
(strictly) admissible solution
of either (\ref{hesskF}), (\ref{rats}), or (\ref{sigmakt}) in
$\Omega = B(O,r_0) \setminus \{O\}$, with
$f(x) \geq c_0 > 0$ near $O$.  The goal is to show that $u$ can be extended to a
H\"older continuous function $u^{\ast} \in C^3(\Omega) \cap
C^{\gamma}(B(O,r_0))$.

We first observe that the Ricci curvature of $g_u = e^{-2u}g$ is strictly positive:
\begin{align*}
Ric(g_u) \geq c_1 g_u,
\end{align*}
where $c_1 = c_1(\min f) > 0$.  This follows from \cite{GVW}; see
Lemma 4.1 and the remark thereafter of \cite{GVRic} for a proof.

Now, in view of Theorem \ref{thm1'} $(i)$, if
\begin{align} \label{gb}
\liminf_{x \rightarrow 0} u(x) > - \infty
\end{align}
then we are done.  Therefore, we must have
\begin{align*}
\liminf_{x \rightarrow 0} u(x) = - \infty,
\end{align*}
and from part $(ii)$ of Theorem \ref{thm1'} it follows that
\begin{align}  \label{blow}
u(x) \leq 2 \log d(x) + C.
\end{align}
Note that since $f$ is not identically zero, the scale invariant
gradient estimate is not necessarily satisfied, so
we only get the upper inequality.
To summarize the idea of the proof, we first show that
(\ref{blow}) implies $g_u$ has geodesics of arbitrary length.
However, since $g_u$ has strictly positive Ricci curvature, this
will yield a contradiction, and it will follow that (\ref{gb}) must hold.

To analyze the behavior of $g_u$ near the singularity let $\{ x^i \}$ denote
conformally flat coordinates centered at $O$, so that $g_{ij} =
\delta_{ij}$ in $\Omega = B(O, r_0) \setminus \{ O \}$.
Let us write the metric as follows:
\begin{align}
g_u = e^{-2u} g = e^{-2 \Psi} \cdot {|x|}^{-4} g =  e^{-2 \Psi} g_{
\star}
\end{align}
where
\begin{align} \label{Psidef}
\Psi =  u - 2 \log {|x|},
\end{align}
\begin{align} \label{gstardef}
g_{\star} = {|x|}^{-4} g.
\end{align}

Changing coordinates to $z^j = |x|^{-2} x^j$, which are
defined on $\mathbf{R}^n \setminus B(0, r_0)$, it is easy to see
that
\begin{align*}
(g_{\star})_{ij}(z) &= \delta_{ij}.
\end{align*}
Therefore,
\begin{align}
\label{decay1} g_u = e^{ -2 \Psi} \delta_{ij} \geq e^{-2 C}
\delta_{ij}.
\end{align}

 The inverted $z$-coordinates are only defined on the complement of a
 ball, so
let us extend $\Psi$ arbitrarily to a function defined on all of
$\mathbf{R}^n$, and consider the metric $\tilde{g} = e^{-2\Psi}
\delta_{ij}$.

\begin{lemma}
The metric $(\mathbf{R}^n, \tilde{g})$ is geodesically complete.
\end{lemma}
\begin{proof}
Let $x_0 \in \mathbf{R}^n$, and let $\zeta(t)$ be a unit-speed
geodesic with $\zeta(0) = x_0$. Assume the maximal domain of
definition of $\zeta$ is $[0,T)$; we want to show that $T =
\infty$. We will make use of the following property of maximal
geodesics: $\zeta: [0,T) \rightarrow \mathbf{R}^n$ must leave every
compact subset of $\mathbf{R}^n$ as $t \rightarrow T$. That is,
given any compact subset $K \subset \mathbf{R}^n$, there exists a
$t_K$ such that $\zeta(t) \in \mathbf{R}^n \setminus K$ for all $t
> t_K$. For a simple proof of this fact, see \cite[page
109]{Petersen}.

Therefore, without loss of generality, we may assume there is a time
$0 < a < T$ such that $\zeta(t) \in \mathbf{R}^n \setminus B(0, R)$
for $t \geq a$. For $t \geq a$, by (\ref{decay1})
\begin{align}
\begin{split}
\label{estC2}
| \dot{\zeta}(t)|^2_{g_w} & = (g_w)_{ij}(\zeta(t)) \dot{\zeta}(t)^i \dot{\zeta}(t)^j \\
& \geq  ( e^{- 2C_1} \delta_{ij}) \dot{\zeta}(t)^i \dot{\zeta}(t)^j\\
& \geq  e^{- 2C_1} | \dot{\zeta}(t)|^2_{0}.
\end{split}
\end{align}

Now let $b \in (a,T)$, so $\zeta(b) \in \mathbf{R}^n \setminus
B(0,R)$.  Since $\zeta$ has unit speed, the length of
$\zeta([0,b])$ is given by
\begin{align} \label{L0b}
\begin{split}
b = L( \zeta([0,b]) = \int_0^b | \dot{\zeta}(t)|_{g_w} dt
&= \int_0^{a}  | \dot{\zeta}(t)|_{g_w} dt + \int_a^b  | \dot{\zeta}(t)|_{g_w} dt \\
&\geq a + e^{-C_2}\int_a^b | \dot{\zeta}(t)|_0 dt.
\end{split}
\end{align}
Since segments minimize distance in the Euclidean metric, we have
\begin{align*}
\int_a^b | \dot{\zeta}(t)|_0 dt \geq |\zeta(b) - \zeta(a)|_0.
\end{align*}
Therefore,
\begin{align}
\label{goody} b - a & \geq e^{-C_1}| \zeta(b) - \zeta(a)|_0.
\end{align}

Now, recall that given any compact set $K \subset \mathbf{R}^n$,
there must be a time $t_K$ with $\zeta(t) \in \mathbf{R}^n
\setminus K$ for $t
> t_K$. Therefore, by choosing a large enough compact set we can
arrange so that $\zeta(b) \in \mathbf{R}^n \setminus K$ and
$|\zeta(b) - \zeta(a)|_0$ is as large as we like.  By
(\ref{goody}), this means we can choose $b$ as large as we like,
i.e., $T = \infty$. It follows that $(\mathbf{R}^n, \tilde{g})$ is geodesically complete.
\end{proof}

 To finish the proof, $\tilde{g}$ is a complete
$C^3$ metric on $\mathbf{R}^n$ which has strictly positive Ricci
curvature outside of a compact set. By Myers' Theorem (see
\cite{Petersen}), this is impossible.  More precisely, take any
constant $N
> 0$, and choose a point $y$ in $\mathbf{R}^n$ with $|y|_0 > N + R$.
By the Hopf Rinow Theorem (which is valid for $C^3$ metrics), there
exists a unit speed minimizing geodesic $\zeta(t) : [0, d(x_0, y)]
\rightarrow \mathbf{R}^n$ with respect to the metric $\tilde{g}$
with $\zeta(0) = 0$. Choose the smallest time $a$ so that
$\zeta(t) \in \mathbf{R}^n \setminus B(0,R)$ for $t \geq a$. From
(\ref{goody}), we have
\begin{align}
d(x_0, y) - a \geq e^{-C_1} | y - R|_0 \geq e^{-C_1} N.
\end{align}
We can therefore find a minimizing geodesic in $\mathbf{R}^n
\setminus B(0,R)$ with arbitrarily long length. But from our
assumption, together with Newton's inequality, and the estimate
$({\bf iv})$, $\tilde{g}$ has strictly positive Ricci curvature $Ric
> c g > 0$, on $\mathbf{R}^n \setminus B(0,R)$, so Myer's Theorem
gives a upper bound on the length (depending only upon $c$), which
is a contradiction.

\begin{remark}
We remark why this proof does not work in the general (non-locally
conformally flat) case, even if we assume $u \in C^4$. The main
point is that since $f$ is not identically zero, a
scale-invariant $C^2$-estimate does not necessarily
hold (see Section \ref{Pointwise} below).
In the course of the proof we used Theorem \ref{thm1'}, which requires
local conformal flatness, while
Theorem \ref{thm1} in the general case requires a scale-invariant $C^2$-estimate
to hold at the singularity.
Also,  Theorem \ref{thm1'} only requires the solution to
be $C^{1,1}_{loc}$, but we must assume $C^3$ to be able to apply
the tools from Riemannian geometry that were used above.
\end{remark}

\section{Integral and H\"older estimates for admissible functions}
\label{inthol}

In this section we prove Theorem \ref{holderthm}. Let $u \in
C^{1,1}_{loc}$ satisfy $A_u \in \overline{\Gamma}_{\delta}$ almost
everywhere in $B(O,r_0),$ and once again denote
\begin{align*}
v = \ds e^{\beta u},
\end{align*}
where $\beta$ is defined in (\ref{betadef}). We first prove
(\ref{optab}), from which the H\"older estimate (\ref{hest2first})
follows.  Afterwards we will consider the conformally flat case.

\begin{proposition}
\label{claim3} Let
\begin{align} \label{p0def}
p_0 = \ds 2 + \frac{1}{\delta} > n.
\end{align}
Then given any $q  < p_0 - 1$, there is a constant $C =
C(\delta,q,n,\mu,g)$ such that
\begin{align} \label{newintgrad2}
\ds \| \nabla v \|_{L^{\frac{n}{(n-1)}q}(B(x_0,r/2))} \leq C(1 +
r^{-1}) \| \nabla v \|_{L^{q}(B(x_0,3r/4))} + C \| v
\|_{L^{q}(B(x_0,r))}.
\end{align}
\end{proposition}

\begin{proof}  We first show that (\ref{hessvcone}) implies that $v$ satisfies a
certain inequality with respect to the $p$-laplacian:

\begin{lemma}
\label{divform2} Let $p_0$ be given by (\ref{p0def}). Then
\begin{align} \label{divpos2}
\nabla_i \big( |\nabla v|^{p_0-2} \nabla_i v \big) \geq - C|\nabla
v|^{p_0-2}
\end{align}
a.e., where $C = C(\delta,n,g)$.
\end{lemma}

\begin{proof}
First, note that (\ref{hessvcone}) implies
\begin{align} \label{vcon2}
\nabla^2 v + \delta \Delta v\ g + \mu v g \geq 0 \quad \mbox{a.e.}
\end{align}
for some constant $\mu = \mu(\delta,n,g) > 0$. Now, a simple
calculation gives
\begin{align*}
\nabla_i \big( |\nabla v|^{p_0-2}  \nabla_i v \big) &= (p_0-2)
\nabla^2 v(\nabla v,\nabla v)|\nabla v|^{p_0-4}  + \Delta v
|\nabla v|^{p_0-2}  \\
&= (p_0-2) \big [ \nabla^2 v + \frac{1}{(p_0-2)}\Delta v\ g + \mu v
g \big ](\nabla v,\nabla v)|\nabla v|^{p_0-4}  \\ & \quad \quad -
\mu(p_0-2) |\nabla v|^{p_0-2}v.
\end{align*}
Since $p_0$ satisfies
\begin{align*}
\ds \frac{1}{(p_0-2)} = \delta,
\end{align*}
from (\ref{vcon2}) we conclude
\begin{align*}
\nabla_i \big( |\nabla v|^{p_0-2} \nabla_i v \big) \geq  - C|\nabla
v|^{p_0-2} v.
\end{align*}
\end{proof}

Now, suppose $2 < p < p_0$.  Then
\begin{align} \label{divID} \begin{split}
\ds |\nabla v|^{2-p} \nabla_i \big( |\nabla v|^{p-2} \nabla v_i
\big) = \frac{p-2}{p_0 - 2} |\nabla v|^{2-p_0} \nabla_i \big(
|\nabla v|^{p_0-2} \nabla v_i \big) + \frac{p_0 - p}{p_0 - 2} \Delta
v.
\end{split}
\end{align}
By Lemma \ref{divform2}, this implies
\begin{align} \label{divIE}
\ds |\nabla v|^{2-p} \nabla_i \big( |\nabla v|^{p-2} \nabla v_i
\big) \geq \frac{p_0 - p}{p_0 - 2} \Delta v - Cv,
\end{align}
which we rewrite as
\begin{align} \label{lapIE}
\ds  \nabla_i \big( |\nabla v|^{p-2} \nabla v_i \big) \geq \frac{p_0
- p}{p_0 - 2} |\nabla v|^{2-p} \Delta v - C|\nabla v|^{2-p} v.
\end{align}

\begin{lemma} \label{deltohess} Suppose $v \in C^{1,1}_{loc}(B(x_0,r))$ satisfies
(\ref{vcon2}).  Then
\begin{align} \label{s2sub}
\Delta v \geq c_0 |\nabla^2 v| - C_1 v
\end{align}
for positive constants $c_0, C_1$ depending on $\delta, n$, and $
g$.
\end{lemma}

\begin{proof}  Choose any point $P \in B(x_0,r)$ at which (\ref{vcon2})
holds.  Let $\{ \nu_1, \nu_2, \dots \nu_n \}$ denote the eigenvalues
of $\nabla^2 v(P)$.  Then $\Delta v(P) = \nu_1 + \cdots + \nu_n$,
and (\ref{vcon2}) implies
\begin{align} \label{vconloc}
\nu_i + \delta \Delta v(P) + \mu v(P) \geq 0
\end{align}
for each $1 \leq i \leq n$.  Summing this inequality for $i \neq j$
\begin{align} \label{vconloc2}
\sum_{i \neq j} \nu_i +  (n-1) \delta \Delta v(P) + \mu (n-1) v(P)
\geq 0,
\end{align}
which gives
\begin{align} \label{vconloc3}
-\nu_j + \big[ 1 + (n-1)\delta \big]\Delta v(P) + \mu (n-1) v(P)
\geq 0.
\end{align}
From (\ref{vconloc}) and (\ref{vconloc3}) we conclude that each
eigenvalue satisfies
\begin{align}
 \big[ 1 + (n-1)\delta \big]\Delta v(P) + C v(P) \geq \nu_j
\geq - \delta \Delta v(P) - C v(P).
\end{align}
Inequality (\ref{s2sub}) follows immediately.
\end{proof}

By (\ref{lapIE}) and (\ref{s2sub}),
\begin{align} \label{hessIE}
\ds  \nabla_i \big( |\nabla v|^{p-2} \nabla v_i \big) \geq
c_0^{\prime} |\nabla v|^{2-p} |\nabla^2 v| - C|\nabla v|^{2-p} v
\end{align}
for positive constants $c_0^{\prime}, C$.  Let $\eta \in
C_0^{\infty}$ be a smooth, non-negative cut-off function supported
in $B = B(x_0,r)$, with $\eta(x) = 1$ in $B(x_0,r/2)$, $\eta(x) = 0$
on $B \setminus B(x_0,\frac{2}{3}r)$, and $|\nabla \eta| \leq
Cr^{-1}$. Multiplying both sides of (\ref{hessIE}) by $\eta$,
integrating and applying the divergence theorem gives
\begin{align*}
\int \eta |\nabla^2 v||\nabla v|^{p-2} \leq C\int |\nabla
\eta||\nabla v|^{p-1} + C \int \eta |\nabla v|^{p-2}v.
\end{align*}
By H\"older's inequality,
\begin{align*}
\int \eta |\nabla v|^{p-2}v \leq \Big( \int \eta |\nabla v|^{p-1}
\Big)^{(p-2)/(p-1)} \Big( \int \eta v^{p-1} \Big)^{1/(p-1)}.
\end{align*}
Therefore, by the properties of $\eta$,
\begin{align} \label{moser}
\int_{B(x_0,r/2)} |\nabla^2 v||\nabla v|^{p-2} \leq C\big(1 +
r^{-1}\big) \int_{B(x_0,3r/4)} |\nabla v|^{p-1} + C \int_{B(x_0,r)}
v^{p-1}.
\end{align}
By the Sobolev imbedding theorem,
\begin{align} \label{sob1} \begin{split}
\Big( \int_{B(x_0,r/2)} |\nabla v|^{\frac{n}{n-1}(p-1)}
\Big)^{(n-1)/n} &\leq C \int_{B(x_0,r/2)} |\nabla |\nabla v|^{p-1}|
+ Cr^{-1} \int_{B(x_0,r/2)} |\nabla v|^{p-1} \\
&\leq C \int_{B(x_0,r/2)} |\nabla^2 v||\nabla v|^{p-2} + Cr^{-1}
\int_{B(x_0,r/2)} |\nabla v|^{p-1} \\
&\leq C\big(1 + r^{-1}\big) \int_{B(x_0,3r/4)} |\nabla v|^{p-1} + C
\int_{B(x_0,r)} v^{p-1}.
\end{split}
\end{align}
Taking $q = p - 1$, this completes the proof of Proposition
\ref{claim3}.
\end{proof}

\begin{corollary} \label{holdercor1}
We have
\begin{align} \label{opt11}
\int_{B(x_0,r/2)} |\nabla v|^p \leq C \big( \int_{B(x_0,r)} v^2
\big)^{p/2}
\end{align}
for any $p < p_{\delta} = \frac{n}{n-1}(p_0 - 1).$
\end{corollary}
\begin{proof}
By tracing (\ref{vcon2}) we see that $v$ satisfies the linear
elliptic inequality
\begin{align*}
\Delta v \geq -Cv
\end{align*}
in $B(x_0,r)$.  Therefore,
\begin{align} \label{pie}
\int_{B(x_0,3r/4)} |\nabla v|^2 \leq C\int_{B(x_0,r)} v^2.
\end{align}
Taking $p = 3 < p_0$ in (\ref{sob1}) and using (\ref{pie}) gives
\begin{align*}
\Big( \int_{B(x_0,r/2)} |\nabla v|^{\frac{2n}{(n-1)}}
\Big)^{(n-1)/n} \leq C \int_{B(x_0,r)} v^2.
\end{align*}
If we now take $p = 1 + 2n/(n-1)$ in (\ref{sob1}) and continue
iterating, we obtain the bound
\begin{align} \label{opt2}
\int_{B(x_0,r/2)} |\nabla v|^p \leq C \big( \int_{B(x_0,r)} v^2
\big)^{p/2}
\end{align}
for any $p < \frac{n}{n-1}(p_0 - 1).$  One can easily check that
\begin{align*}
p_{\delta} = \frac{n}{n-1}(p_0 - 1) = \frac{ n(1 + \delta)}{(n-1)
\delta} > n.
\end{align*}
\end{proof}

The estimate (\ref{optab}) follows from (\ref{opt2}) by
interpolation.  This gives the H\"older estimate (\ref{hest2first}):
\begin{align*}
 \| v \|_{C^{\alpha}(B(x_0,r/2))} \leq C
\int_{B(x_0,r)} |v|,
\end{align*}
for any $\alpha < \gamma = \frac{1 + (2-n)\delta}{1 + \delta}$.
\vskip.1in

To complete the proof of Theorem \ref{holderthm}, it only remains to
show why in the $LCF$ case we can take $\alpha = \gamma$ in
(\ref{hest2first}).  However, this follows from the proof of Theorem
\ref{thm1}: recall that in conformally flat coordinates $\{ x^i \}$
inequality \ref{hessvcone}) is equivalent to
\begin{align} \label{localcone1}
\left( \partial_i \partial_j v + \beta v A_{ij} \right) \in
\overline{\Gamma}_{\delta}
\end{align}
a.e.  Defining $W(x) = v(x) + \Lambda |x|^2$ as before, with
$\Lambda >> 0$ sufficiently large, we can verify that $W$ is a
$\delta$-convex.  After convolving, we obtain a smooth
$\delta$-convex mollification $W_h$, and by Theorem \ref{unifhold}
$W_h$ satisfies the estimate (\ref{hest2first}) with $\alpha =
\gamma$.  Letting $h \rightarrow 0$, we obtain the same estimate for
$W$ and $v$.

\begin{remark}
 In the locally conformally flat case, we note the proof of (\ref{optab}) only
requires the $\delta$-convexity of $v$, so the estimate
(\ref{w1pnew}) follows.
\end{remark}

\bibliography{IMRN_references}
\end{document}